\documentclass[pdflatex,sn-basic,iicol]{sn-jnl}

\usepackage{graphicx}%
\usepackage{multirow}%
\usepackage{amsmath,amssymb,amsfonts}%
\usepackage{amsthm}%
\usepackage{mathrsfs}%
\usepackage[title]{appendix}%
\usepackage{xcolor}%
\usepackage{textcomp}%
\usepackage{manyfoot}%
\usepackage{booktabs}%
\usepackage{algorithm}%
\usepackage{algorithmicx}%
\usepackage{algpseudocode}%
\usepackage{listings}%

\usepackage{subeqnarray}
\usepackage{empheq}
\usepackage{mathtools}
            
\usepackage{float}
\usepackage{subfig}
\usepackage{anyfontsize}

\usepackage{csml}

\usepackage{lineno}

\renewcommand{\u}{\bm{u}}

\newcommand{\rvec}{\bm{r}}

\newcommand{\M}{\bm{M}}
\newcommand{\A}{\bm{A}}

\DeclareMathOperator*{\minimize}{minimize}

\begin{document}

\title{Robust topology optimisation of lattice structures with spatially correlated uncertainties}

\author*[1,2]{\fnm{Ismael} \sur{Ben-Yelun}}\email{i.binsenser@upm.es}
\author[1]{\fnm{Ahmet Oguzhan} \sur{Yuksel}}\email{aoy21@cam.ac.uk}
\author[1]{\fnm{Fehmi} \sur{Cirak}}\email{fc286@cam.ac.uk}

\affil[1]{\orgdiv{Department of Engineering}, \orgname{University of Cambridge}, \orgaddress{Trumpington Street}, \city{Cambridge}, \postcode{CB2 1PZ}, \country{UK}}

\affil*[2]{\orgdiv{E.T.S. de Ingeniería Aeronáutica y del Espacio}, \orgname{Universidad Politécnica de Madrid}, \orgaddress{Pza. Cardenal Cisneros 3}, \city{Madrid}, \postcode{28040}, \country{Spain}}

\abstract{The uncertainties in material and other properties of structures are often spatially correlated. We introduce an efficient technique for representing and processing spatially correlated random fields in robust topology optimisation of lattice structures. Robust optimisation takes into account the statistics of the structural response to obtain a design whose performance is less sensitive to the specific realisation of the random field. We represent Gaussian random fields on lattices by leveraging the established link between random fields and stochastic partial differential equations (SPDEs). The precision matrix, i.e. the inverse of the covariance matrix, of a random field with Mat\'ern covariance is equal to the finite element stiffness matrix of a possibly fractional PDE with a second-order elliptic operator. We consider the finite element discretisation of the PDE on the lattice to obtain a random field which, by design, takes into account its geometry and connectivity. The so-obtained random field can be interpreted as a physics-informed prior by the hypothesis that the elliptic PDE models the physical processes occurring during manufacturing, like heat and mass diffusion. Although the proposed approach is general, we demonstrate its application to lattices modelled as pin-jointed trusses with uncertainties in member Young's moduli. We consider as a cost function the weighted sum of the expectation and standard deviation of the structural compliance. To compute the expectation and standard deviation and their gradients with respect to member cross-sections we use a first-order Taylor series approximation. The cost function and its gradient are computed using only sparse matrix operations. We demonstrate the efficiency of the proposed approach using several lattice examples with isotropic, anisotropic and non-stationary random fields and up to eighty thousand random and optimisation variables.
}

\keywords{robust optimisation, random fields, lattice structures, topology optimisation}

\maketitle


\clearpage

%
\section{Introduction \label{sec:intro}}
%

%
\subsection{Motivation}
%
The inherent random variations in material properties and geometry resulting from manufacturing have adverse implications for the performance and design of optimised structures in general and, more specifically, lattice structures considered in this paper. Lattice structures in the form of micro-architected lattices, or mechanical metamaterials, have recently attracted much attention due to their superior mechanical performance compared to conventional bulk solids~\citep{fleck2010micro, gibson1999cellular,zheng2014ultralight}. Technological advances in additive manufacturing, especially over the past decade, have made it possible to produce lattice structures of various topologies with geometric features spanning multiple length-scales~\citep{gibson2021additive,schaedler2016architected,maconachie2019slm}. The vast design space afforded by additive manufacturing can be utilised using structural optimisation, in particular by optimising the member cross-sections and topology.   

In structural optimisation it is essential to consider the uncertainties in material properties \citep{doltsinis2004robust,asadpoure2011robust,da2017stress}, geometry \citep{guest2008structural,chen2011new,jansen2015robust}, loading \citep{ben1997robust,kogiso2008robust,zhao2014robust,torres2021robust} and manufacturing parameters \citep{wang2011projection,schevenels2011robust}. A design obtained using conventional deterministic optimisation is usually sensitive to variations in material properties, loading and geometry, often eradicating the optimality of the structure. In robust optimisation, uncertainties are taken into account by considering a cost function consisting of the weighted sum of the expectation and standard deviation of standard cost functions used in deterministic optimisation; see~\cite{park2006robust,beyer2007robust,schueller2008computational} for an overview. Consequently, robust optimisation yields a structure which is optimal in terms of a chosen cost function and is less sensitive to the uncertain properties of the actual structure. For completeness, we note that alternative approaches for considering uncertainties, including reliability-based structural optimisation, are not discussed further in this paper. 
\begin{figure*}
	\centering
	\subfloat[Original design]{\includegraphics[width=0.33\linewidth]{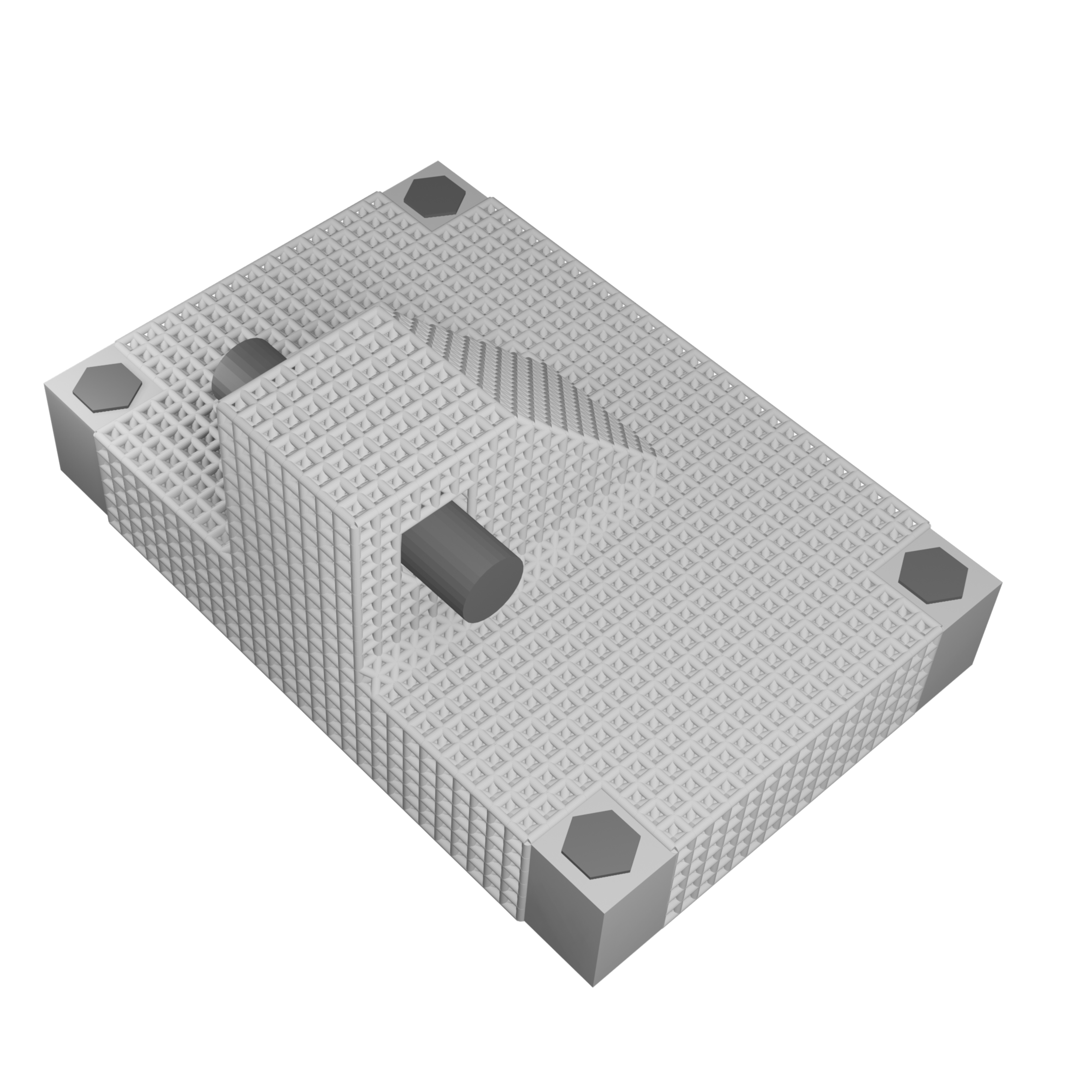}\label{intro_gbc}}
	\hfill
	\subfloat[Random member Young's moduli]{\includegraphics[width=0.33\linewidth]{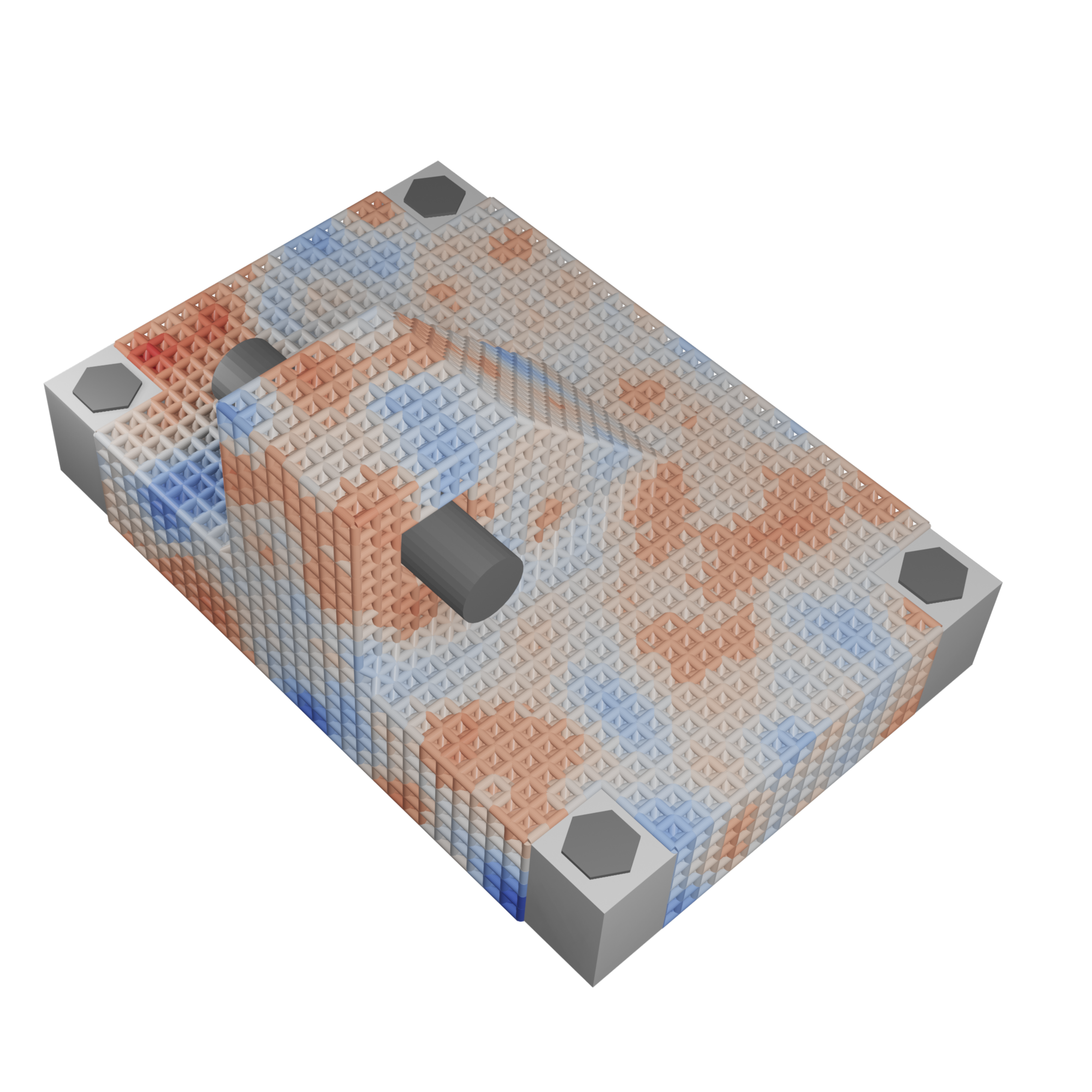}\label{intro_sample}}
	\hfill
	\subfloat[Optimised design]{\includegraphics[width=0.33\linewidth]{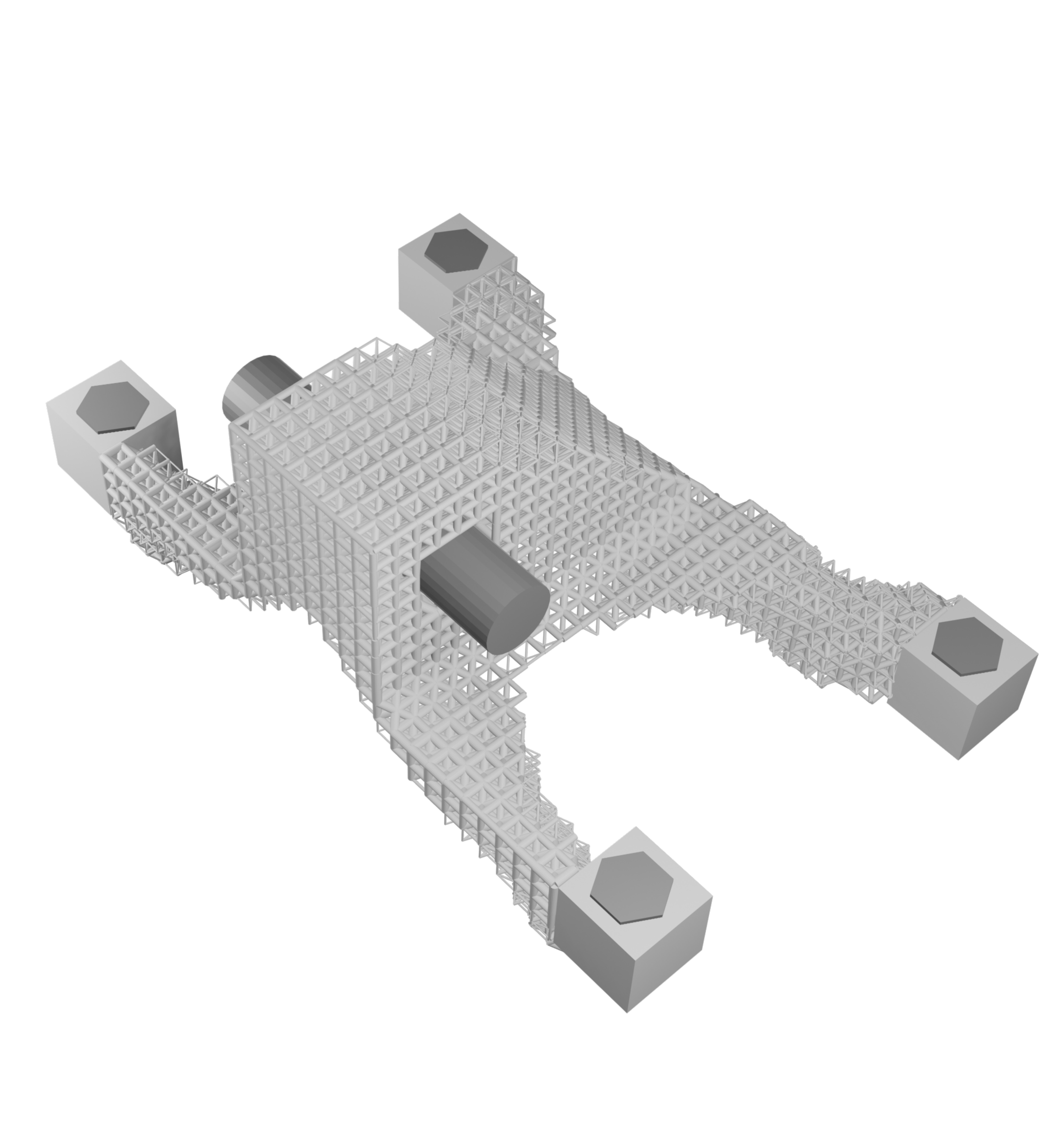}\label{intro_opt}}
	\caption{Robust topology optimisation of a bracket consisting of body-centred cubic unit cells. The bracket is fixed at the four corners of the base plate and an external load is applied at the horizontal shaft attached to the vertical plate. For further details on this example see Section \ref{sec:bracket}. }
	\label{fig:1_summary}
\end{figure*}

%
\subsection{Related research}
%
The variations in material properties and geometry of the lattice structure can be represented as spatial random fields, or processes. We assume in this paper that the respective random fields are approximated as Gaussian fields with a mean and a covariance. The covariance, or kernel, function describes the correlation between the field values at different positions in space, and consequently depends on the specific manufacturing and assembly process used. Especially, in the case of the Young's modulus the covariance function must depend on the manufacturing process and the partial differential equations corresponding to the underlying physical processes. For instance, in additive manufacturing heat is used to melt/sinter powders or to melt filaments to make three-dimensional lattices layer by layer~\citep{ gibson2021additive}.  Consequently, it is expedient to consider the Laplace operator modelling the heat flow in establishing the covariance function of random fields. Indeed the formal link between PDEs and covariance properties of random fields is widely known in spatial statistics~\citep{whittle1954stationary, rue2005gaussian,lindgren2011explicit,lindgren2022spde}. For instance, \cite{ whittle1954stationary} showed that a Gaussian random field with the widely used Mat\'ern covariance function is equivalent to the solution of a stochastic PDE with a Gaussian white noise as the forcing. The respective PDE with a second-order elliptic operator is of fractional order. In other words, after discretising it with finite elements, the fractional power of the corresponding stiffness matrix is equal to the precision matrix, \ADDED{i.e. the inverse of the covariance matrix,} of the random field~\citep{lindgren2011explicit}.  The so-obtained covariance matrix can be interpreted as a physics-informed prior describing the spatial correlation in the random field. 

The benefits of the PDE representation of the covariance include that the resulting matrices are usually sparse and that the PDE can be easily generalised to anisotropic and non-stationary (non-homogenous) random fields on bounded and non-Euclidean domains. For lattices, the correlated random field can be obtained by considering the discretisation of the PDE operator on the lattice. That is, the respective global diffusion matrix is assembled from the one-dimensional element diffusion matrices and the global mass matrix from the one-dimensional element mass matrices. Subsequently, the covariance matrix is determined using standard procedures developed for the non-lattice case~\citep{bolin2020rational,koh2023stochastic}. The sketched approach is somewhat similar to the graph kernels recently proposed in statistical machine learning~\citep{borovitskiy2020matern,borovitskiy2021matern,nikitin2022non, bolin2023regularity}. 

Robust optimisation is computationally significantly more costly than deterministic optimisation because of the dependence of cost and, possibly, constraint functions on the second-order statistics of the finite element solution. The probability density of the solution is the push-forward of the probability density of the material and other random parameters by the lattice solution operator. Hence, even though the input probability densities are Gaussians, the probability density of the finite element solution is usually non-Gaussian and is analytically intractable.  
A straightforward but computationally prohibitive approach is to use Monte Carlo (MC) sampling \citep{haldar2000probability} to estimate the expectation and variance of the objective and cost functions and their gradients (sensitivities) in each optimisation step. MC sampling alone is hardly scalable beyond a handful of random variables and is consequently of limited practical use. A multitude of approaches has been explored to reduce the number of design variables and the number of required finite element solutions to estimate the statistics of the cost and constraint functions by MC sampling. Perhaps the most common approach is to reduce the dimensionality of random fields by using a stochastic series expansion, including Karhunen-Lo\'eve ~\citep{sudret2000stochastic} or the expansion optimal linear estimator (EOLE)  \citep{li1993optimal}. In structural optimisation series expansions have been used, e.g., in~\cite{schevenels2011robust,lazarov2012topology,lazarov2012atopology,jansen2013robust,zhao2014robust}. 

Alternatively, the number of finite element evaluations can be drastically reduced with a truncated Taylor series expansion of the finite element solution with respect to the random parameters to estimate expectation and variance of cost and constraint functions~\citep{doltsinis2004robust, guest2008structural,asadpoure2011robust, lazarov2012topology}.  Usually, only the first-order linear term in the Taylor series expansion is considered so that the probability density of the structural response is a Gaussian, given that the input densities are Gaussians, yielding the so-called perturbation methods \citep{kriegesmann2019robust}. \ADDED{The perturbation methods inevitably lead to truncation errors. The coefficient of variation of the random parameters must  be small to obtain meaningful approximations.}

Instead of the series expansion, other approaches for solving stochastic PDEs~\citep{ matthies2005galerkin}, like the stochastic Galerkin method using a polynomial chaos expansion (PCE)~\citep{ghanem1991stochastic},  can be utilised to estimate the second order statistics of the random variables and fields required in robust optimisation~\citep{tootkaboni2012topology,keshavarzzadeh2016gradient}. Although the PCE based approaches are appealing, they are unsuitable for problems with a large number of random variables.

For completeness, methods which replace the finite element model with a pre-trained surrogate model, or a metamodel, are worth mentioning as well~\citep{ jin2003use, santner2003design, sobester2008engineering}. The surrogate model can be, for instance, a Gaussian process regression model~\citep{williams2006gaussian} or a more sophisticated probabilistic machine learning model~\citep{ vadeboncoeur2022deep} trained on a dataset consisting of finite element solutions. The surrogate model is very fast to evaluate and consequently makes MC sampling practical. More recent techniques to speed up the MC sampling include the use of stochastic gradient-based methods prominent in machine learning~\citep{de2020topology}.  

%
\subsection{Contributions}
%
In this paper, we introduce a novel approach for the robust topology optimisation of lattice structures with up to several tens of thousands of members, see Figure~\ref{fig:1_summary}.  We focus on correlated random fields, especially their efficient representation and consideration within gradient-based robust optimisation. The lattice structure is modelled as a pin-jointed truss and consists of a large number of unit cells each of which have a small number of members. This approximation is sufficient for technologically relevant stretch-dominated lattices with appealing mechanical properties; see, e.g. \cite{deshpande2001foam}. We optimise the topology of the lattice by taking the cross-sectional areas of the members as design variables~\citep{xiao2022infill, christensen2008introduction}. Although the proposed approach is general, we demonstrate its application to lattices with uncertainties in member Young's moduli, and the cost function is the weighted sum of the expectation and standard deviation of the compliance. 

The probability density of the member Young's moduli is given by a multivariate Gaussian distribution with a covariance matrix which we determine by generalising the SPDE representation of random fields to  lattices. Because Young's moduli are associated with the members of the lattice, we consider the adjoint lattice for discretising the stochastic PDE. By interpreting the lattice structure as a graph, its adjoint (or, line) graph is another graph representing the adjacencies between the members of the lattice~\citep{ harary1960some}. The adjoint lattice is defined by placing the vertices of the adjoint graph at the centroids of the members of the original lattice.  As mentioned, the precision matrix of the random field is equal to the stiffness matrix of the discretised fractional PDE with a second-order elliptic operator on the adjoint lattice. The precision matrix is sparse when the PDE is of integer order and can be approximated as sparse when non-integer order~\citep{bolin2020rational,koh2023stochastic}. We consider only integer order operators in this paper. 

Furthermore, we use a first-order Taylor series approximation to determine the second-order statistics of the random compliance.  The mean and standard deviation of the compliance and their gradients with respect to the cross-sectional areas are  determined using only sparse matrix operations given that the precision matrix of the random Young's moduli is sparse. Consequently, the proposed approach can be applied to lattices with a very larger number of members and joints. 

In addition, we apply a density filter to avoid scattered members (similar to checkerboards in continuum structures) and  penalisation to obtain cross-sectional areas within a certain prescribed range~\citep{xiao2022infill}. 

\section{Robust design optimisation \label{sec:robustDesign}}
%
In this section, we review the robust topology optimisation of linear elastic truss structures. That is, the minimisation of the member cross-sectional areas for a prescribed loading and material volume.
The Young's moduli of the members is assumed to be random and spatially correlated. 

%
\subsection{Problem formulation \label{sec:formulation}}
%
The equilibrium equation for a pin-jointed truss structure with~$n_e$ members and~$n_d$ degrees of freedom is given by 
\begin{equation} \label{eq:equilibrium}
	\vec K\left(\vec s, \vec r\right)\u\left(\vec s, \vec r\right) = \vec f  \, , 
\end{equation}
where~$\vec f \in \mathbb R^{n_d}$ and~$\vec u \in \mathbb R^{n_d}$ are the displacement and force vectors, and~$\vec K \in \mathbb R^{n_d \times n_d}$ is the positive definite symmetric global stiffness matrix. The stiffness matrix and the displacements depend on member design variables and Young's moduli collected in the vectors~$\vec s \in \mathbb R^{n_e}$ and~$\vec r  \in \mathbb R^{n_e}$, respectively. The design vector~$\vec s$ consists of the relative (pseudo) density of the members defined as  
\begin{equation}
	s_e = \frac{A_e - A_{\min}}{A_{\max}- A_{\min}} \, . 
\end{equation}
Here,~$A_e$ is the cross-sectional area of the member~$e$ and~$A_{\max}$ is the maximum allowable cross-sectional area. Furthermore,~$A_{\min}$ is a small algorithmic parameter  to avoid the ill-conditioning of the stiffness matrix  when~\mbox{$A_e \approx 0$}. The cross-sectional area and Young's modulus of each member is uniform across its length. 

The Young's modulus of the members is random and has the multivariate Gaussian density
\begin{equation} \label{eq:densityYoungsModuli}
	 p(\vec r) = \mathcal{N}\left(\overline{\vec r},  \vec{C}_{r} \right ) \, ,
\end{equation}
where~$\overline{\vec r} \in  \mathbb R^{n_e}$ is its expected value and~$\vec C_r \in  \mathbb R^{n_e \times n_e}$  its positive definite covariance matrix. Both are defined as
\begin{equation} \label{eq:densityYoungsModuliMeanAndCov}
	\overline{\vec r} = \expect [\vec r] \, , \quad \vec C_r = \expect \left [ (\vec r- \overline{\vec r})(\vec {r}- \overline{\vec {r}})^{\trans}\right ] \, ,
\end{equation}
using the expectation operator~$\expect [\cdot]$. The covariance matrix is usually dense and becomes a diagonal matrix when the components of the random vector are uncorrelated, and in case of Gaussian densities statistically independent.  Furthermore, the precision matrix~$\vec Q_r$ is defined as the inverse of the covariance matrix, i.e.~\mbox{$\vec Q_r = \vec C_r^{-1}$.} 
  
The robust compliance optimisation problem for the truss structure with random member Young's moduli can now be stated as 
 \begin{subequations}
\begin{align}
\displaystyle\minimize_{\vec s} & \quad F(\vec s) = \alpha_E\overline{J}(\vec s) + \alpha_{\sigma} \sigma_J (\vec s) \label{eq:optimisation_subeq_cost}\\
 \mathrm{subject\;to} & \quad \vec K\left(\vec s, \vec r \right)\vec u \left(\vec s,\vec r\right) = \vec f \, , \label{eq:optimisation_subeq_equilibrium}\\
& \quad g(\vec s) = V(\vec s) - V_{\max} \leq 0 \, ,  \label{eq:optimisation_subeq_constraint} \\
& \quad \mathbf{0} \leq \vec s \leq \vec I_{n_e} \,   . 
\label{eq:optimisation_subeq_designspace}
\end{align}
\label{Eq:robust_optimisation}%
\end{subequations} 
The compliance is defined as 
\begin{equation}
\label{eq:compliance}
\begin{split}
J(\vec s, \vec r) &= \vec f \cdot \vec u(\vec s, \vec r) \\
&= \vec u(\vec s, \vec r) \cdot \vec K(\vec s, \vec r) \vec u(\vec s, \vec r) \, ,
\end{split}
\end{equation}
where\ADDED{the dot denotes inner product}and the last expression is obtained by considering the total potential of the structure at equilibrium. The cost function~$F(\vec s)$ is the sum of the expectation and variance,  
\begin{subequations}
	\label{eq:exactCostMeanCov}
	\begin{align}
	\overline{J}(\vec s) &= \expect \left[J(\vec s, \vec r)\right] \,  , \label{eq:exactCostMeanCov:a} \\
	\sigma_J(\vec s) &= \sqrt{\expect \left[ \left (J(\vec s, \vec r) - \overline{J}(\vec s) \right )^2\right]} \, , \label{eq:exactCostMeanCov:b}
	\end{align}
\end{subequations}
weighted by the two user-chosen parameters~$\alpha_E, \, \alpha_{\sigma} \in \mathbb R^+$. The choice of the two parameters will be detailed in Section~\ref{sec:examples} on examples. The constraint~\eqref{eq:optimisation_subeq_constraint} ensures that the volume~\mbox{$V(\vec s) = \vec {l} \cdot \vec {a}(\vec s)$} of the optimised structure is below the prescribed volume~$V_{\text{max}}$, where $\vec {l} \in \mathbb R^{n_e}$ and~$\vec {a}(\vec s) \in \mathbb R^{n_e}$  are the vectors of member lengths and  cross-sectional areas, respectively.  The last constraint~\eqref{eq:optimisation_subeq_designspace} enforces that the relative density~$s_e$ of the members is between~$0$ and~$1$.

%
\subsection{Perturbation approximation \label{sec:perturbation}}
%
The fast evaluation of the cost function and its gradient are critical in iterative structural optimisation. As specified in~\eqref{eq:optimisation_subeq_cost} the cost function is the weighted sum of the expectation and standard deviation of the compliance. We use a first-order Taylor series expansion of the displacements to approximate both. 

The first-order series expansion of the random displacement vector~$\vec u (\vec s, \vec r)$ around the expected value~$\overline{\vec r}$ of the random vector is given by
\begin{equation} \label{eq:solExpand}
\u(\vec s, \vec r) \approx \u(\vec s, \overline{\vec r}) + \dfrac{\partial\u(\vec s,\overline{\vec r})}{\partial\vec r}\left(\vec r- \overline{\vec r}\right) \, .
\end{equation}
The gradient of the solution herein is derived by differentiating the equilibrium equation~\eqref{eq:equilibrium} yielding
\begin{equation}
\begin{split}
\begin{aligned}
\label{Eq:equilibrium_derivative}
&\frac{\partial}{\partial \vec r} \left ( \vec K(\vec s, \vec r)\vec u (\vec s,\vec r) - \vec f \right ) = \vec 0 \quad \Rightarrow \\
&\frac{\partial\vec u(\vec s, \overline{\vec r})}{\partial \vec r} = -\vec K(\vec s, \overline{\vec r})^{-1}   \frac{\partial \vec K(\vec s, \overline{\vec r})}{\partial\vec r}\vec u (\vec s, \overline{\vec r}) \,  .
\end{aligned}
\end{split}
\end{equation}
We use in the following the same symbols for the approximate vectors and matrices like their exact counterparts to avoid a proliferation of symbols. 

Owing to the linear transformation property of Gaussian vectors, see, e.g. \cite{williams2006gaussian}, the (approximate) probability density of the displacement vector is following~\eqref{eq:solExpand},~\eqref{Eq:equilibrium_derivative} and~\eqref{eq:densityYoungsModuli} also a Gaussian,  i.e.,   
\begin{equation}
	 p\left(\vec u(\vec s) \right) = \mathcal{N}\left(\overline{\vec u}(\vec s), \vec {C}_{u} (\vec s)\right) \, . 
\end{equation}
According to the series expansion~\eqref{eq:solExpand}, we obtain for the expected value
\begin{equation}
\begin{split}
\overline{\vec{u}}(\vec s) &= \expect\left[\vec u\left(\vec s, \overline{\vec r}\right) + \dfrac{\partial\vec{u}\left(\vec s, \overline{\vec r}\right)}{\partial\vec{r}}\left(\vec r-\overline{\vec r}\right)\right] \\
&= \vec u\left(\vec s, \overline{\vec r}\right) \, ,
\end{split}
\end{equation}
and for the covariance matrix
\begin{align}\label{eq:solCovariance}
&\vec {C}_{u} (\vec s) = \expect \left[(\vec u(\vec s, \vec r)  - \overline{\u} (\vec s))(\vec u(\vec s, \vec r) - \overline{\u} (\vec s))^\trans \right] \nonumber\\
& = \dfrac{\partial \u(\vec s, \overline{\vec r}) }{\partial \vec r}\vec C_r \left(\dfrac{\partial \vec u(\vec s, \overline{\vec r})}{\partial \vec r}\right)^{\trans}  \, .
\end{align}
 Furthermore, using one more time the linear transformation property of Gaussian vectors, we arrive at the probability density of the compliance 
 \begin{subequations} \label{eq:approxCost}
\begin{align}
 	& p(J) = \mathcal{N}\left(\overline{J} (\vec s),  \sigma_J (\vec s)^2\right) \,  ,  
	\intertext{with the mean and variance}
	& \overline{J} (\vec s) = \vec f \cdot \overline{\u} (\vec s) \, , \label{eq:approxCostMean} \\
    & \sigma_J(\vec s)^2 = \vec f \cdot \vec {C}_{\u}(\vec s, \overline{\vec r}) \vec f \, . \label{eq:approxCostVar}	
\end{align}
\end{subequations}
For later derivations, we can can express the variance more compactly as
\begin{equation} \label{eq:approxCostVarCompact}
		\sigma_J(\vec s)^2 = \frac{\partial J (\vec s, \overline{\vec r})}{\partial \vec r } \cdot \vec {C}_{ r}  \frac{\partial J (\vec s, \overline{\vec r})}{\partial \vec r } \, , 
\end{equation}
where it has been taken into account that

\begin{equation}
    \dfrac{\partial J(\vec s, \overline{\vec r})}{\partial \vec r}  = -\vec u(\vec s, \overline{\vec r})\cdot \dfrac{\vec K(\vec s, \overline{\vec r})}{\partial \vec r}\vec u(\vec s, \overline{\vec r}) \, .
\end{equation}
%
\subsection{Sensitivities \label{sec:sensitivities}}
%
We determine the gradient, i.e. sensitivity, of the cost function~\eqref{eq:optimisation_subeq_cost} using the first-order approximation of the random displacements~\eqref{eq:solExpand}.  To this end, it is beneficial to consider the cost function as the sum of the member contributions, i.e.,   
\begin{equation} \label{Eq:objective_function_sensitivity}
		F (\vec s)  = \alpha_E \sum_{e=1}^{n_e}  \overline{J}_e (\vec s) + \alpha_{\sigma} \sum_{e=1}^{n_e} \sigma_{J,e} (\vec s) \, . 
\end{equation}
By recourse to the approximate expectation~\eqref{eq:approxCostMean} and the compliance~\eqref{eq:compliance}, the contribution of member~$e$ to the gradient of the expected compliance is given by
\begin{equation}
	 \dfrac{\partial \overline{J}_e(\vec s)}{\partial s_e}  = -\u_e (\vec s, \overline{\vec r}) \cdot \dfrac{\partial \vec K_e (s_e, \overline{\vec r}) }{\partial s_e}\u_e  (\vec s, \overline{\vec r}) \, .
\end{equation}
Here,~$s_e$ denotes the design variable,~$\vec K_e$ the element stiffness matrix and~$\vec u_e$ the nodal displacement vector of member~$e$. Similarly, we can obtain from the approximate variance~\eqref{eq:approxCostVarCompact} the member~$e$'s contribution to the gradient of the standard deviation as
\begin{equation}
\begin{aligned}
\dfrac{\partial\sigma_{J,e} (\vec s)}{\partial s_e} &= \dfrac{1}{2\sigma_{J,e} (\vec s) }\dfrac{\partial \left (\sigma_{J,e} (\vec s) \right )^2} {\partial s_e} \\
&= \dfrac{1}{ \sigma_{J,e} (\vec s)}  \dfrac{\partial^2J (\vec s, \overline {\vec r})}{\partial\vec r\partial s_e} \cdot \vec {C}_{r} \dfrac{\partial J (\vec s, \overline{\vec r})}{\partial\vec r}  ,
\end{aligned}
\end{equation}
where
\begin{align}
&\dfrac{\partial^2 J (\vec s, \overline{\vec r})}{\partial\vec r\partial s_e} = -2\u(\vec s, \overline{\vec r}) \cdot \dfrac{\partial \vec K (\vec s, \overline{\vec r})}{\partial\vec r}\dfrac{\partial\u (\vec s, \overline{\vec r})}{\partial s_e} \nonumber \\
&\quad\quad -\u_e(\vec s, \overline{\vec r}) \cdot \dfrac{\partial^2\vec K_e (s_e, \overline{\vec r})}{\partial\vec r\partial s_e}\u_e (\vec s, \overline{\vec r}) \, .  
\end{align}
The required gradient of the displacement vector is determined by solving for each~$s_e$ the linear equation system 
\begin{equation}\label{eq:displacement_derivative}
\vec K(\vec s, \overline{\vec r}) \dfrac{\partial\u (\vec s, \overline{\vec r})}{\partial s_e} = 
-\dfrac{\partial \vec K(\vec s, \overline{\vec r})}{\partial s_e} \vec u (\vec s, \overline{\vec r}) \, .
\end{equation}

Finally, for gradient-based optimisation also the gradient of the volume constraint~\eqref{eq:optimisation_subeq_constraint} is needed, which can be straightforwardly computed as
\begin{equation} \label{Eq:constraint_sensitivity}
\dfrac{\partial g (\vec s)}{\partial s_e} = l_e \dfrac{\partial A_e (s_e)}{\partial s_e} \, ,
\end{equation}
where~$l_e$ and~$A_e$ are the length and cross-sectional area of member~$e$.

%
\subsection{Penalisation and filtering \label{Sec:filter_penalisation}}
%
We use in some of the examples included in this paper a penalisation approach akin to the SIMP method~\citep{bendsoe1999material} to obtain optimised structures with a desired cross-sectional area distribution. Specifically, we aim to reduce the number of members with a relative density~$s_e < s^*$. The  prescribed parameter~$s^*$ will usually depend on the specific manufacturing process used, and herein is set to $s^* = 0.5$. Figure~\ref{fig:spline_plots} shows the used the piecewise $\mathcal{C}^1$  continuous penalised density~$\tilde s_e (s_e)$ obtained by combining a univariate B-spline curve and a line. The B-spline is of polynomial degree four, has six control points, and the begin and end points are interpolating.\ADDED{The choice of different penalisation functions has been explored by~\cite{xiao2022infill}.}
\begin{figure}
	\includegraphics[width=\linewidth]{spline}\label{spline}
	\caption{Penalisation function mapping the relative density~$s_e$ to the penalised relative density~$\tilde s_e$. The piecewise penalisation function consists of a B-spline curve and a line, which are smoothly connected at~$s_e^*=0.5$. The dashed polygon and its corners represent the control polygon of the B-spline curve. The penalisation ensures that densities~$s_e =0$ and~$s_e > s^*$ are preferred. 	\label{fig:spline_plots}}
\end{figure}

Furthermore, we regularise the optimisation problem by applying a filter to the design variables to avoid checkerboard-like instabilities known from topology optimisation for continuum structures. The used cone kernel function evaluated at the centroidal coordinates $\vec {c}_i$ and $\vec {c}_j$ of the members $i$ and $j$ reads
\begin{equation}
w_{ij} = \max\left\lbrace0, R- \lVert \vec c_i- \vec c_j \rVert \right\rbrace \, , 
\end{equation}
where~$R$ is the prescribed filter length and~$\| \cdot \|$ the Euclidean distance.  With the obtained filtering weights the filtered design variable of a member~$e$ is given by
\begin{equation} \label{Eq:filter}
\hat{s}_e = \left ( \sum_i \dfrac{w_{ei}s_i}{l_i} \right ) \bigg / \left ( \sum_i\dfrac{w_{ei}}{l_i} \right) \, .
\end{equation}
As proposed in~\cite{xiao2022infill} this definition includes the member length~$l_i$ rendering the cost function gradient a measure of the strain energy density (elastic energy per unit volume) rather than the total compliance. To apply filtering to the gradients, the vector of design variables~$\vec s$ in the cost function~\eqref{Eq:objective_function_sensitivity} and the constraint~\eqref{Eq:constraint_sensitivity} is replaced with its filtered counterpart~$\hat {\vec s}$. Subsequently, their gradients are determined using the chain rule, that is,  
\begin{subequations}
	\begin{align}
	\dfrac{\partial F (\hat{\vec s})}{\partial s_e} &= \sum_k\dfrac{\partial F (\hat{\vec s})}{\partial \hat{s}_k}\dfrac{\partial \hat{s}_k}{\partial s_e} \, , \\
	\dfrac{\partial g(\hat{\vec s})}{\partial s_e} &= \sum_k\dfrac{\partial g (\hat{\vec s})}{\partial \hat{s}_k}\dfrac{\partial \hat{s}_k}{\partial s_e} \, .
	\end{align}
\end{subequations}
Differentiation of the kernel function~\eqref{Eq:filter} yields the derivative
\begin{equation}
\dfrac{\partial \hat{s}_k}{\partial s_e} =\dfrac{w_{ke}}{l_e}   \dfrac{1}{\sum_i\dfrac{w_{ki}}{l_i}}.
\end{equation}

%
\section{Spatial random fields \label{sec:randomFields}}
%
We now turn to the computation of the multivariate Gaussian probability density of the member Young's moduli. After reviewing the stochastic PDE representation of Mat\'ern fields and its finite element discretisation on bounded Euclidean domains, we present its generalisation to lattice structures. Our focus is on obtaining the precision matrix corresponding to a generalised Mat\'ern field with prescribed hyperparameters, such as length-scale, smoothness, variance and anisotropy.  
%
\subsection{Mat\'ern fields on unbounded Euclidean domains \label{sec:maternGeneral}}
\begin{figure}
	\centering
	\subfloat[Mat\'ern covariance function ]{\includegraphics[width=1.06\linewidth]{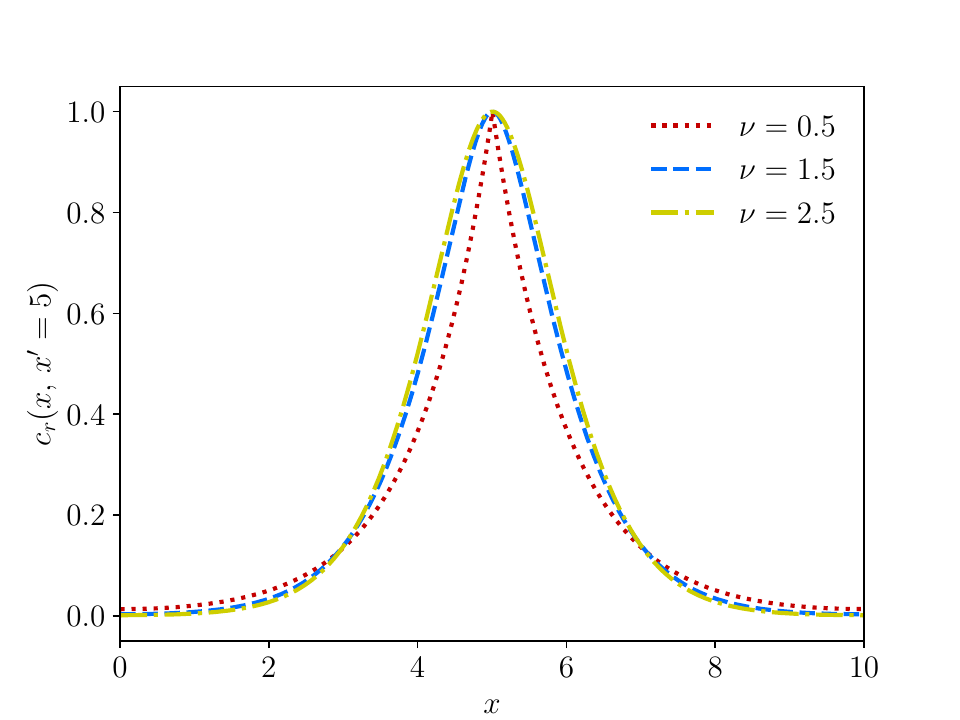}\label{fig:1d_nu}} \\
	\subfloat[Samples drawn from Mat\'ern Gaussian processes]{\includegraphics[width=1.06\linewidth]{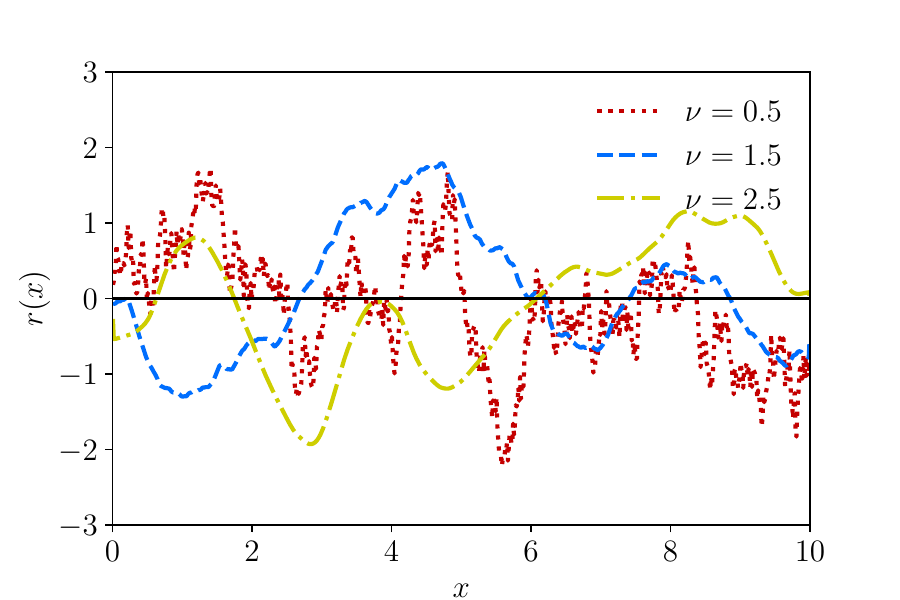}\label{fig:1d_len}} 
	\caption{Mat\'ern covariance function $c_r(x, x' = 5)$ with $\sigma = 1$, $\ell=1$, $\nu \in \{ 0.5,1.5,2.5 \}$, and three samples drawn from respective Mat\'ern Gaussian processes with a zero mean. The black line in (b)  represents the zero mean.}
	\label{Fig:3_matern_covariance_1d}
\end{figure}

%
We consider in~$\mathbb R^d$, $d \in \{1,2,3 \}$, a zero-mean Gaussian process,  or a (spatial) Gaussian random field,
\begin{equation}
	r(\vec  x)\sim \set{GP} \left (0, c_r \left(\vec  x, \vec  x' \right)  \right)  \, ,
\end{equation}
with the Mat\'ern covariance function $c_r(\vec x, \vec x')$ defined as 
\begin{equation}
\begin{aligned}
\label{eq:exact_matern}
&c_r(\vec x, \vec x') =  \mathbb{E}[r(\vec x)r(\vec x')] \\
&= \frac{\sigma^2}{2^{\nu - 1}\Gamma(\nu)}\left(\kappa \lVert \vec x-\vec x'\rVert\right)^{\nu} K_{\nu}\left(\kappa \lVert \vec x-\vec x'\rVert\right) \, ,
\end{aligned}
\end{equation}
with $\sigma\in \mathbb{R}^{+}$ being the standard deviation, $\nu\in\mathbb{R}^{+}$ the smoothness parameter, $\ell\in\mathbb{R}^{+}$ the length-scale parameter, $\Gamma$ the Gamma function and $K_{\nu}$ the modified Bessel function of the second kind of order $\nu$. The parameter $\kappa$ is defined as
\begin{equation}
\kappa = \dfrac{\sqrt{2\nu}}{\ell} \, .
\label{eq:kappa}
\end{equation}
Figure~\ref{Fig:3_matern_covariance_1d} shows the covariance function for smoothness parameters $\nu \in \{ 0.5,1.5,2.5 \}$ and realisations of the corresponding Gaussian zero-mean random fields. As apparent, the covariance function and the respective samples become smoother with increasing $\nu$. 

According to~\cite{whittle1954stationary} and~\cite{lindgren2011explicit}, the random field~$r(\vec x)$ is the solution of the stochastic partial differential equation (SPDE)  
\begin{equation} \label{eq:spde}
\left(\kappa^2 - \Delta\right)^{\beta}r(\vec  x) = \dfrac{1}{\tau}g(\vec  x) \, ,
\end{equation}
where $\Delta$ is the Laplace operator, $g(\vec  x)$ a Gaussian white noise process 
\begin{equation} \label{eq:GPwhiteNoise}
g(\vec  x) \sim \mathcal{GP}(0,\;\delta (\vec  x-\vec  x')) \, ,
\end{equation}
and the remaining two parameters are defined as
\begin{equation}
\beta = \dfrac{\nu}{2} + \dfrac{d}{4}\, ,  \;   \tau^2 = \dfrac{\Gamma (\nu)}{\sigma^2 \Gamma (\nu + d/2) (4\pi)^{d/2}\,\kappa^{2\nu}} \, .\label{subeq:tau}
\end{equation}
Although the exponent~$\beta$ can take any value \mbox{$\beta > d/4$}, we choose in this paper~\mbox{$\beta \in \mathbb Z_{>0}$} to avoid the need to consider fractional derivatives. 

Following~\cite{lindgren2011explicit}, for a smoothness parameter~$\beta > 1$ the SPDE is solved using the recursion 
\begin{subequations} \label{eq:spdeRecurs}
\begin{align}
	\left(\kappa^2 - \Delta\right) r^{(1)}(\vec  x) &= \dfrac{1}{\tau}g(\vec  x) \,  ,  \label{eq:spdeRecurs1} \\  
	\left(\kappa^2 - \Delta\right) r^{(k)}(\vec  x) &=  r^{(k-1)}(\vec  x)  \,   \label{eq:spdeRecurs2}, 
\end{align} 
\end{subequations}
i.e.~$r(\vec x) \equiv r^{(\beta)} (\vec x)$ with $k = 2, \dotsc \beta $. Without going into details, a straightforward finite element discretisation of the PDEs in the recursion yields 
\begin{subequations} \label{eq:spdeFErecurs}
\begin{align}
\left(\kappa^2 \M + \A \right) \vec r^{(1)}  &= \dfrac{1}{\tau}\vec g \, ,   \label{eq:spdeFErecurs1} \\
\left(\kappa^2 \M + \A \right) \vec r^{(k)}  &= \M \vec r^{(k-1)} \, ,
\label{eq:spdeFErecurs2}
\end{align}
\end{subequations}
where~$\M$ is the (lumped) mass matrix,~$\A$ the stiffness matrix and~$\vec g$ a Gaussian white noise vector  with the density 
\begin{equation}
	\vec g \sim \set{N}( \vec 0, \vec M) \, ;
\end{equation}
see~\cite{koh2023stochastic}.

Owing to the linear transformation property of Gaussian random vectors, the solution vector~$\vec r$ is a Gaussian and has the probability density 
\begin{equation} \label{eq:spdeFE}
\vec r  \sim \set {N}(\vec 0, \vec C_r) = \set {N}\left(\vec 0, \vec Q_r^{-1}\right) \, , 
\end{equation}
with the precision matrix 
\begin{equation}\label{eq:precision}
	\vec Q_{r} = \tau^2 \sqrt{ \vec M} \left(\kappa^2 \vec{I} + \sqrt{\vec M^{-1}} \A \sqrt{\vec M^{-1}} \right)^{2\beta} \sqrt{\vec M} \, . 
\end{equation}
\ADDED{We refer to Appendix~\ref{appen:SPDEdiscret}  for further details on this derivation and the computation of the precision matrix. }

%
\subsection{Mat\'ern fields on lattices \label{sec:maternLattices}}
%
The stochastic PDE representation of random fields is instrumental in generalising Mat\'ern fields to lattices. The naive approach of evaluating the covariance function~\eqref{eq:exact_matern} by introducing the coordinates of lattice points does not take into account the topology of the lattice, and simply replacing the distance with the geodesic, or shortest-path distance, does not always lead to a valid, i.e. positive definite, covariance matrix~\citep{gneiting1998simple}. However, starting from the finite element discretised random field~\eqref{eq:spdeFE} and replacing the mass and stiffness matrices with the corresponding matrices of the lattice structure does indeed lead to a positive definite covariance matrix. 

\begin{figure}
    \centering
    \includegraphics[width=0.9\linewidth]{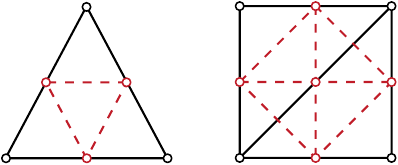}
    \caption{Two lattices  (solid black) and their adjoint lattices (dashed red). The vertices of the adjoint lattice are placed at the centroid of the members of the original lattice.}
    \label{fig:line_graph}
\end{figure}

As noted, Young's modulus of a member is assumed to be constant along its length, so that the members attached to a  joint may have different Young's moduli. This suggests using the adjoint lattice instead of the actual lattice for determining the precision matrix of the random vector~$\vec r$. That is, we interpret the lattice as a graph $\mathcal{G(V, E)}$ where the joints and members of the lattice correspond to the vertices $\mathcal{V}$ and edges $\mathcal{E}$, respectively. In the adjoint graph $\mathcal{L(G)}$, defining the adjoint lattice, the edges of $\mathcal{G}$ become the vertices of~$\mathcal{L(G)}$ \citep{harary1960some}. The adjoint lattice joints are placed at the centroids of the lattice members.  In Figure~\ref{fig:line_graph}, two example lattices and their adjoints are depicted.  

It is straightforward to determine the mass and stiffness matrices of the so-defined adjoint lattice. For instance, for a lattice chain in~$\mathbb R^1$ the  lumped mass and stiffness matrices of a member~$e$ are given by 
\begin{equation} \label{eq:1dMassStiffnessMat}
	\M_e =   \dfrac{h_e}{2}\begin{pmatrix}
	1 & 0 \\ 0 & 1
	\end{pmatrix} \, , \quad
	\A_e =  \dfrac{1}{h_e}\begin{pmatrix*}[r]
	1 & -1 \\ -1 & 1
	\end{pmatrix*} \, , 
\end{equation}
where~$h_e $ is the length of the member. After assembling the element matrices  as in standard finite elements, the obtained matrices are introduced into~\eqref{eq:precision} to obtain the precision matrix for the adjoint lattice.  As usual, the mass and stiffness matrices of a member embedded in~$\mathbb R^2$ or~$\mathbb R^3$ are obtained by applying a coordinate transformation to~\eqref{eq:1dMassStiffnessMat}.

\begin{figure}
	\centering
	\subfloat[Two diagonals per unit cell \label{fig:isotropic_samples_a}]{\includegraphics[width=0.31\textwidth]{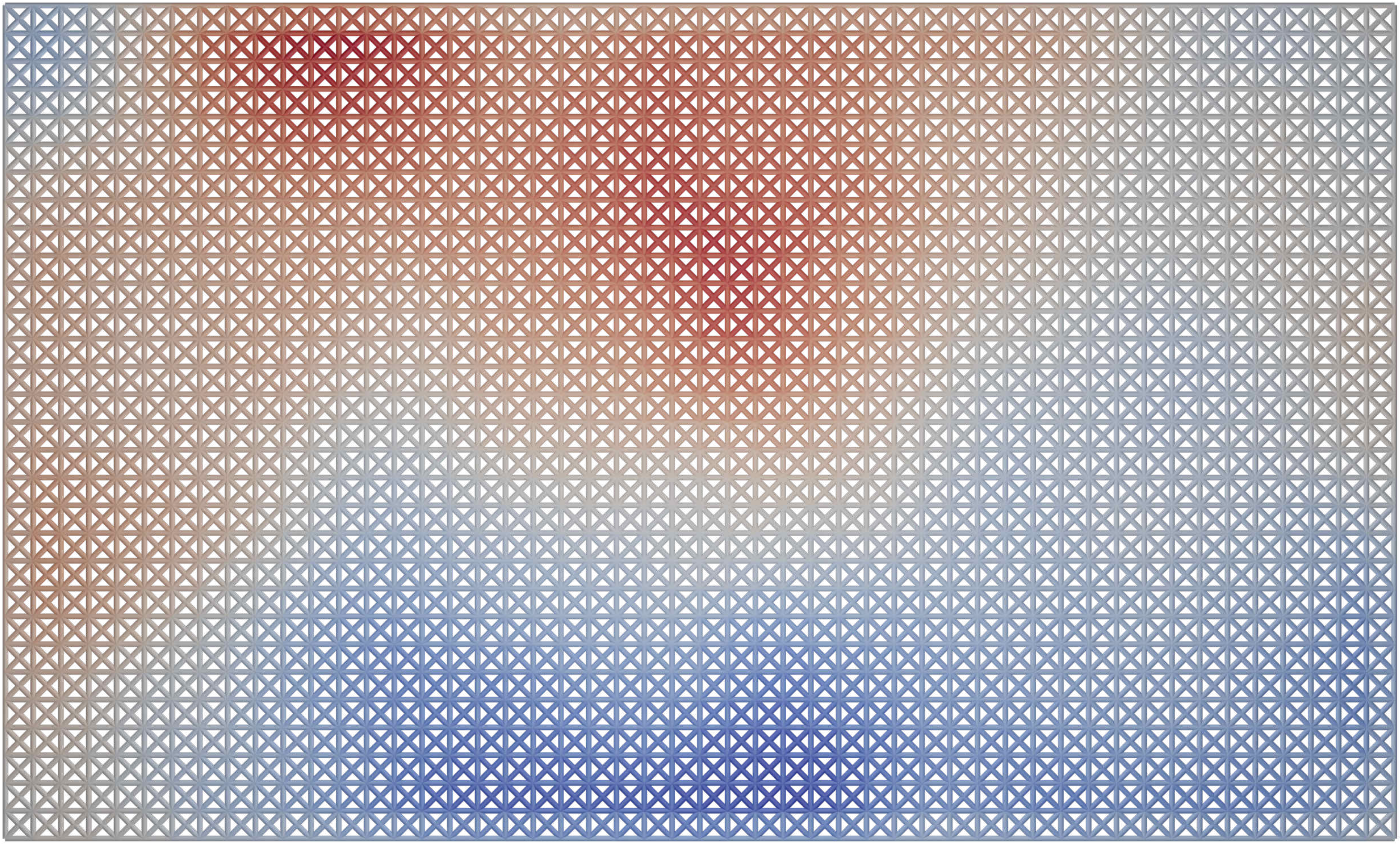}\label{fig:ill_sample_plate}}
	\\
	\subfloat[Two diagonals per unit cell  \label{fig:isotropic_samples_b}]{\includegraphics[width=0.31\textwidth]{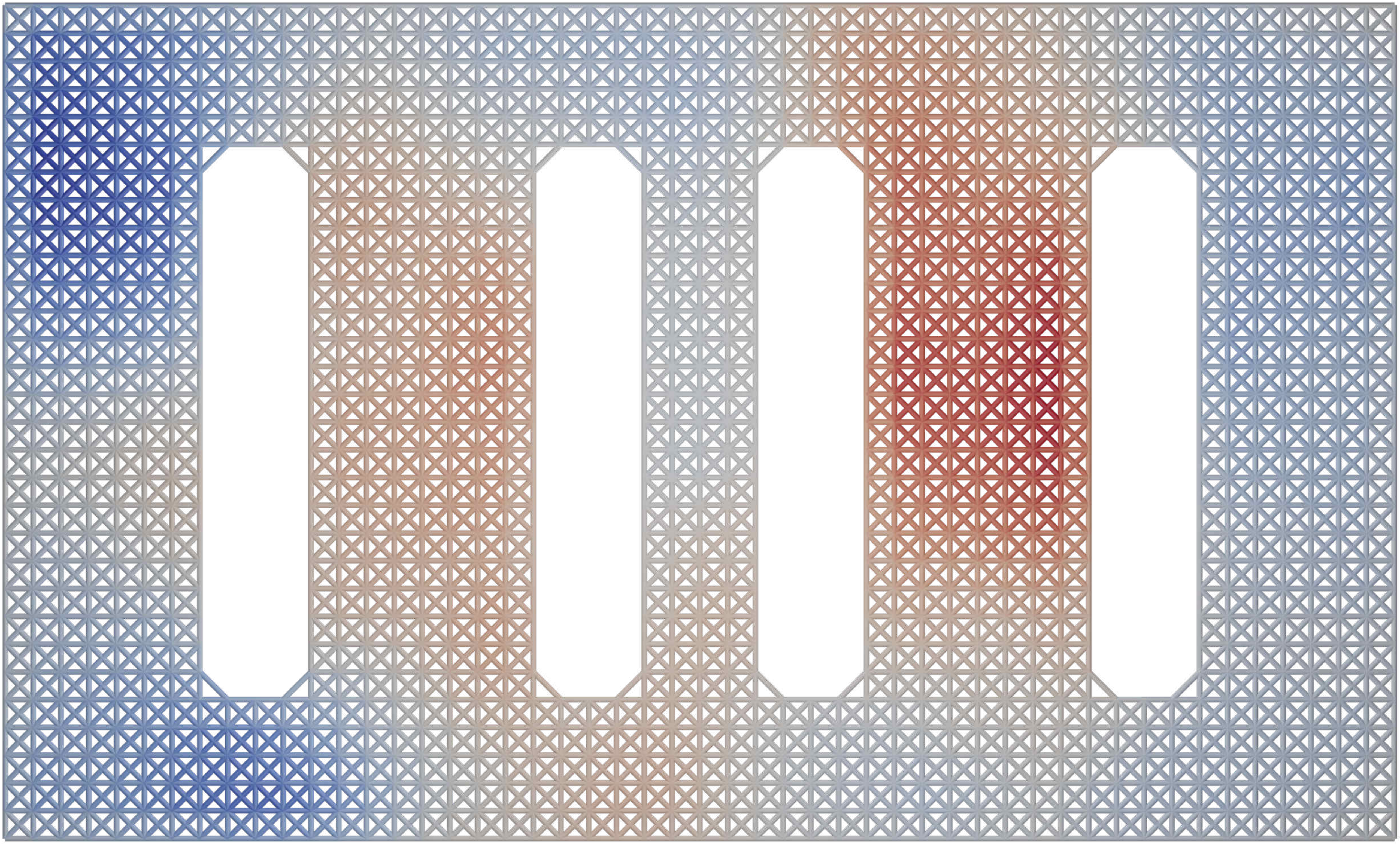}\label{fig:ill_window_2_diag}}
	\\
	\subfloat[One diagonal per unit cell  \label{fig:isotropic_samples_c}]{\includegraphics[width=0.31\textwidth]{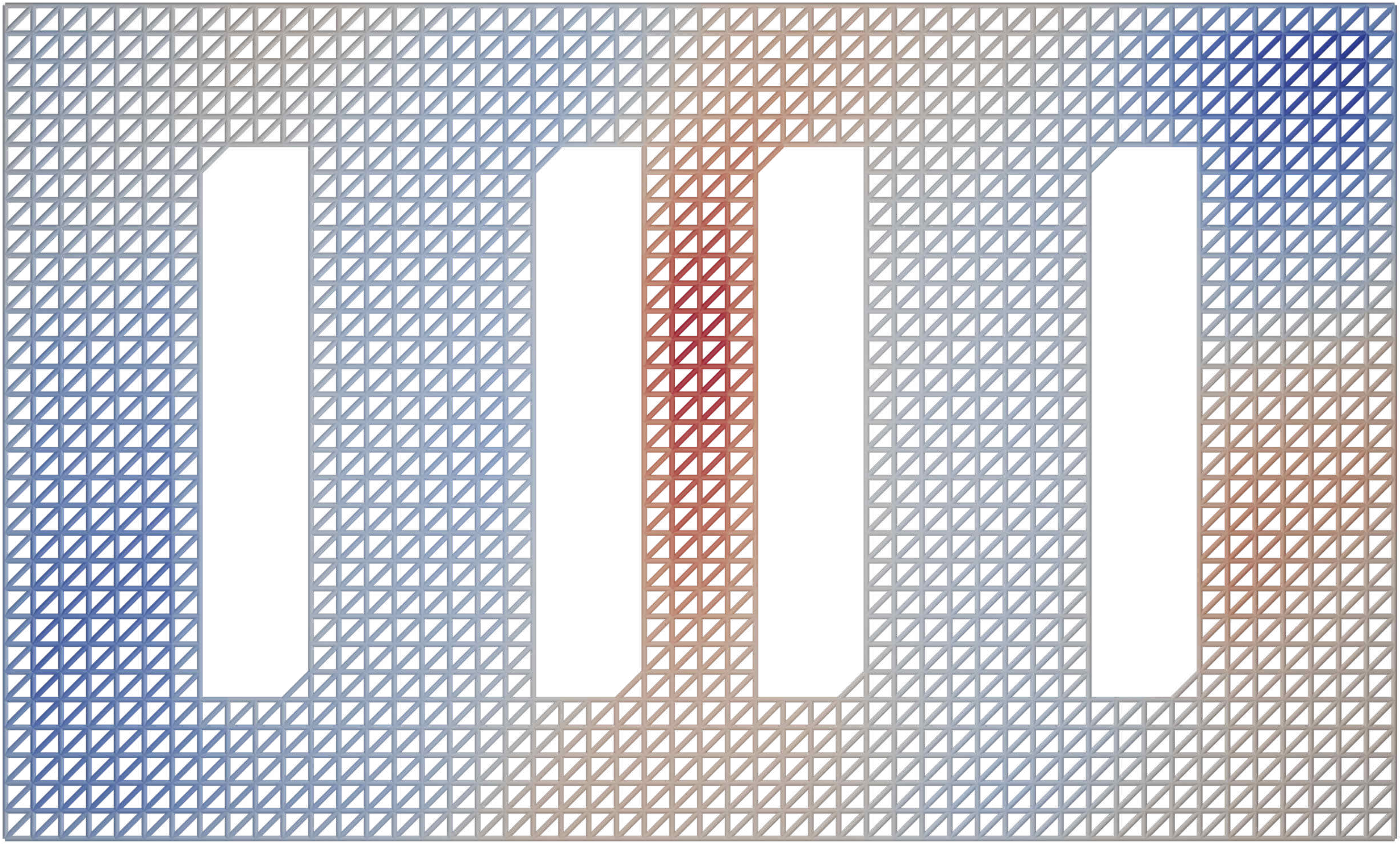}\label{fig:ill_window_1_diag}}
	\caption{Random fields sampled from SPDE-defined Gaussian probability densities. The size of each lattice is~$50 \times 30$ and the size of each unit cell is~$1 \times 1$. The lattice structures in (a) and (b) have unit cells with two diagonals and the lattice structure in (c) has unit cells with a single diagonal.  In all case the length-scale parameter is~$\ell = 20$, the smoothness parameter is~$\nu =2.0$, the standard deviation is~$\sigma = 0.1 \overline{r}$, and the mean is~$\overline { r} = 100$.}
    \label{fig:isotropic_samples}
\end{figure}

Figure~\ref{fig:isotropic_samples} depicts the isocontours of random fields sampled from SPDE-defined Gaussian probability densities. As apparent, the obtained random fields depend on the geometry and topology of the lattice structure and its unit cell connectivity.  The two lattices in Figures~\ref{fig:isotropic_samples_a} and~\ref{fig:isotropic_samples_b} have the same unit cell connectivity with two diagonal members. Comparing the isocontours in Figures~\ref{fig:isotropic_samples_a} and~\ref{fig:isotropic_samples_b},  the four holes in Figure~\ref{fig:isotropic_samples_b} lead to a change in the overall appearance of the random field. Specifically, the random field values at two opposing ends of a hole are very different despite their proximity in terms of Euclidean distance. This is unsurprising given that we determine the random field by solving the system of equations~\eqref{eq:spdeFErecurs} containing the lattice stiffness and mass matrices and a Gaussian white noise forcing. We chose as boundary conditions for the nodes adjacent to holes homogenous Neumann and the nodes at the outer boundary as infinite domain; see~\cite{roininen2014whittle,khristenko2019analysis} for discussions on the choice of boundary conditions in case of Euclidean domains. 

 Furthermore, Figure~\ref{fig:isotropic_samples_c} shows a lattice structure with the same geometry and topology but with each unit cell having only one diagonal. This leads again to a change in the overall appearance of the random field in comparison to the lattice structure with two diagonal members in each cell. The observed dependency of the random field on the geometry, topology and unit cell connectivity agrees with our hypothesis that the SPDE models the diffusion processes, i.e. heat or mass diffusion, during manufacturing. Consequently, the obtained random fields represent a sound physics-informed prior for properties of lattices manufactured, e.g., by additive manufacturing or casting. 
%
\subsection{Non-standard Mat\'ern fields on lattices \label{sec:nonstandard_maternLattices}}
%
The SPDE representation of random  fields makes it possible to easily model non-stationary and anisotropic random fields without altering the introduced overall solution process. Specifically, a non-stationary random field is modelled by choosing a spatially varying length-scale~$\ell(\vec x)$ yielding the SPDE parameters 
\begin{subequations}\label{eq:kappa_tau_nh}
	\begin{align}
	\kappa(\vec{x}) &= \dfrac{\sqrt{2\nu}}{\ell(\vec x)} \, , \\  
	\tau^2(\vec x) &= \dfrac{\Gamma (\nu)}{\sigma^2 \Gamma (\nu + d/2) (4\pi)^{d/2}\,(\kappa(\vec x))^{2\nu}} \, .
	\end{align} 
\end{subequations}

To model anisotropic random fields, e.g. resulting from the particulars of the manufacturing process used, we replace the first equation~\eqref{eq:spdeFErecurs1} in the recursion~\eqref{eq:spdeFErecurs} by   
\begin{equation}\label{eq:aniso_spde}
	\left(\kappa^2 \M + \A \right) \rvec = \dfrac{1}{\tau} \vec D^{-1} \vec g \, .
\end{equation}
Here,~$\vec D$ is a diagonal matrix with the components 
\begin{equation}\label{eq:aniso_diffusion}
	D_{ee} = d_{\vert\vert} \vec n  \cdot \vec t_e  + d_{\perp}(1- \vec n \cdot \vec t_e) \, , 
\end{equation}
where~$\vec t_e \in \mathbb R^{d}$ is the unit tangent to member~$e$ of the lattice structure,~$\vec n \in \mathbb R^d$ is a prescribed unit vector describing the direction of anisotropy, and  $d_{\vert\vert}$ and $d_{\perp} \in \mathbb R^{+} $ are two parameters. In practice, the anisotropy vector~$\vec n$ and the two parameters~$d_{\vert\vert}$ and $d_{\perp}$ will  depend on the specific manufacturing process and protocol used. 

\begin{figure}
	\centering
	\subfloat[$\ell = 3$]{\includegraphics[width=0.31\textwidth]{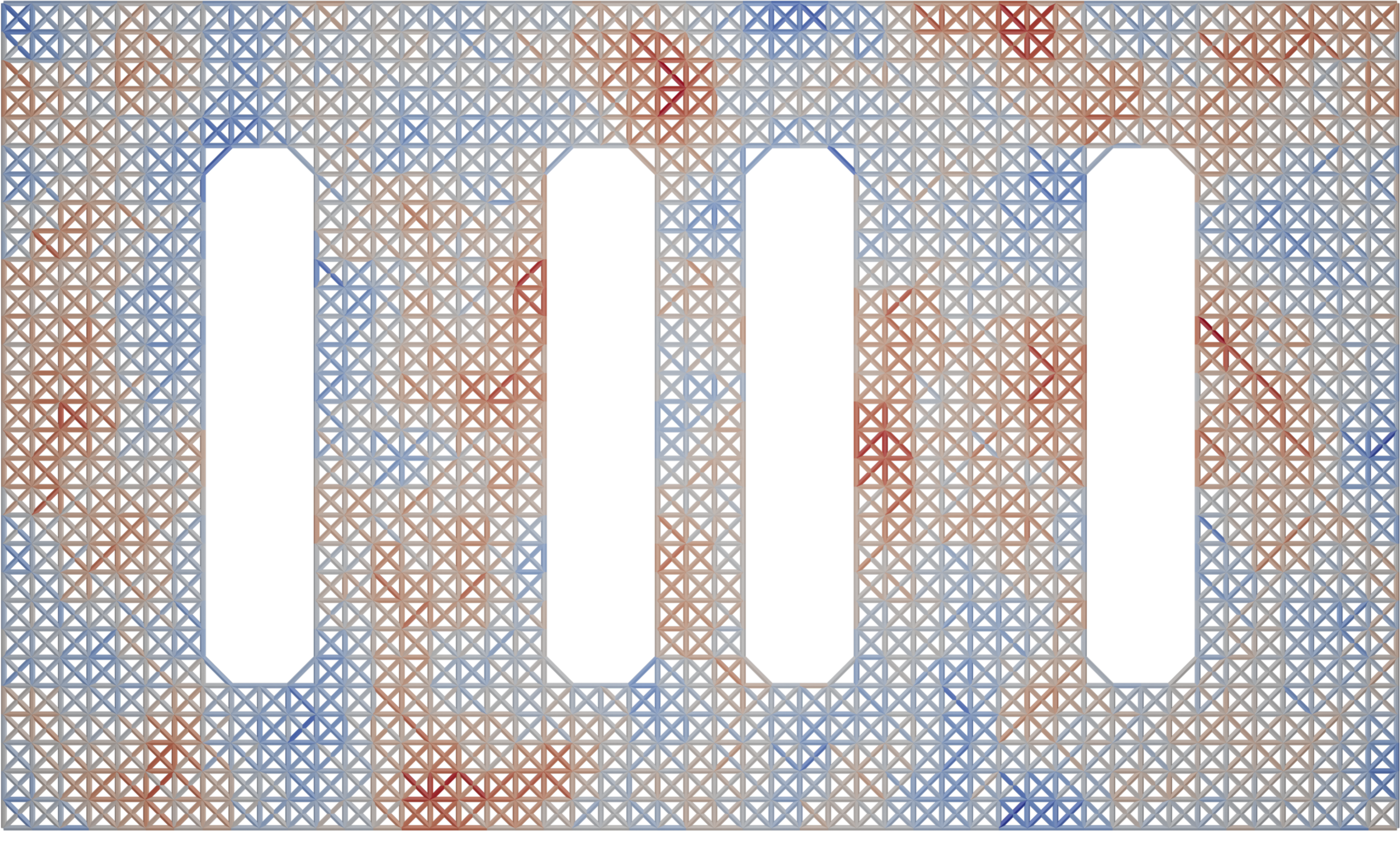}\label{fig:nonh_a}}
	\\
	\subfloat[$\ell(x) = 3+17x/50$]{\includegraphics[width=0.31\textwidth]{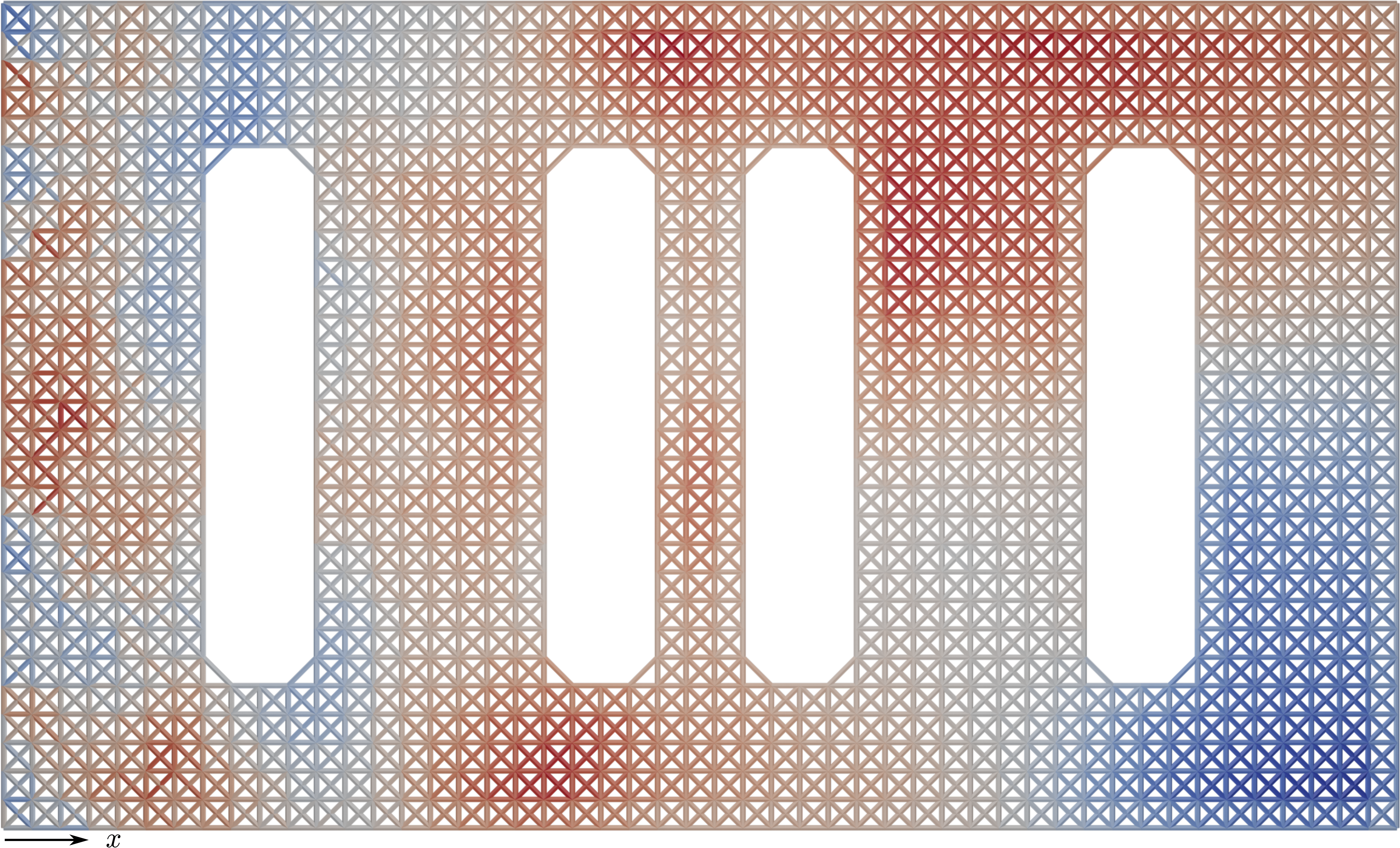}\label{fig:nonh_b}}
	\\
	\subfloat[$\ell = 20$]{\includegraphics[width=0.31\textwidth]{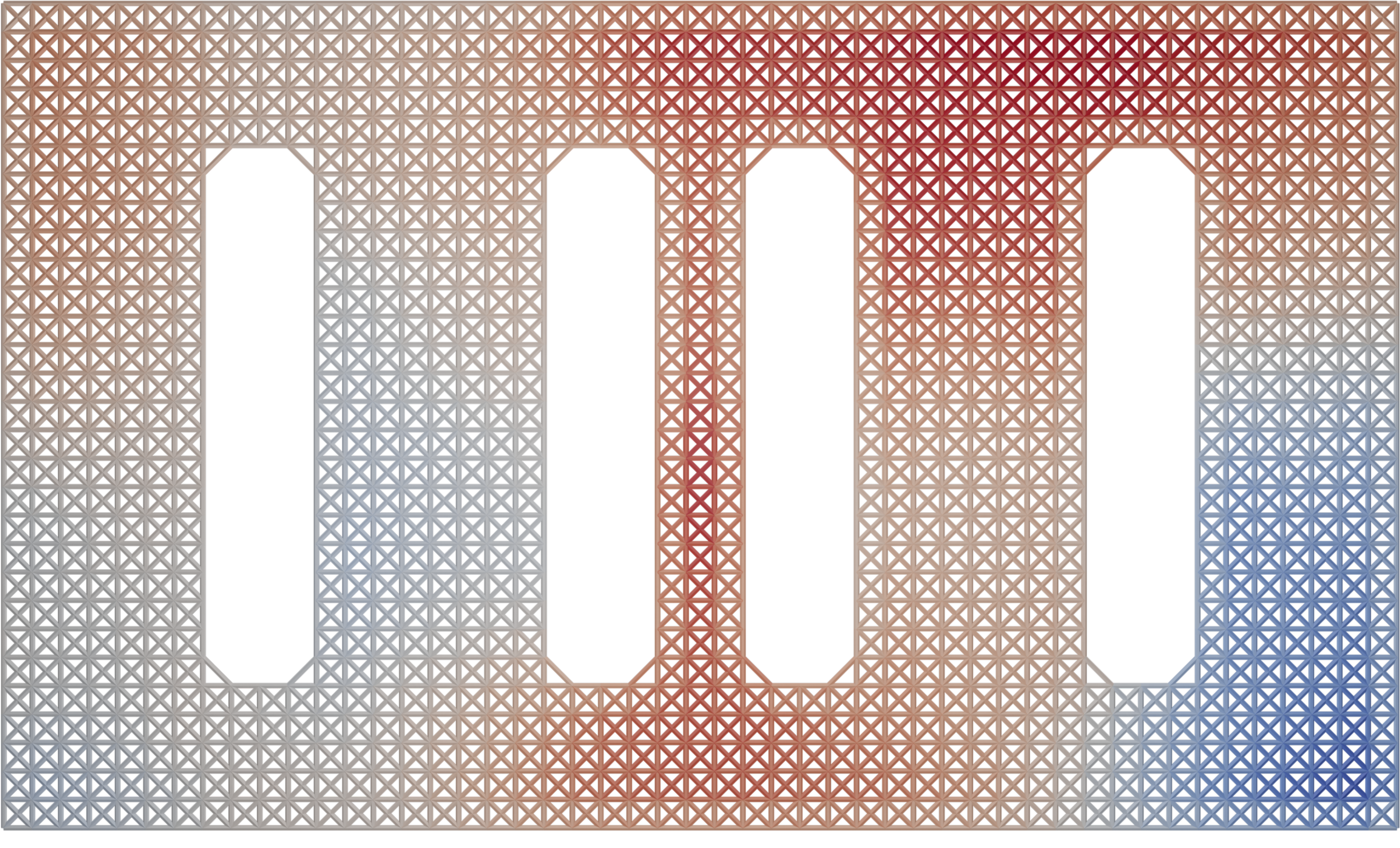}\label{fig:nonh_c}}
	\caption{Stationary and non-stationary random fields sampled from SPDE-defined Gaussian probability densities.   The size of each lattice is~$50 \times 30$, the size of each unit cell is~$1 \times 1$ and each unit cell has two diagonals. In (a) and (c) the length-scale is uniform throughout the domain and in (b)  it increases linearly from~$\ell=3$ at the left boundary to~$\ell=20$ at the right boundary.  The smoothness is~$\nu =2.0$, the standard deviation is~$\sigma = 0.1 \overline{r}$, and the mean is~$\overline {r} = 100$.}
	\label{fig:omega_nh}
\end{figure}

Figure~\ref{fig:omega_nh} shows non-stationary random fields sampled from SPDE-defined Gaussian probability densities. The length-scale of the random field is uniform in Figures~\ref{fig:nonh_a} and~\ref{fig:nonh_c} (\mbox{$\ell=3$} and \mbox{$\ell=20$}, respectively). In Figure~\ref{fig:nonh_b} the non-uniform length-scale increases linearly from~$\ell=3$  at the left end to~$\ell= 20$ at the right end. The increase of the length-scale is easily discernible from the isocontours of the random field. 

%
\section{Examples \label{sec:examples}}
%
%
We proceed to demonstrate the utility and efficiency of the proposed robust optimisation approach with selected examples. In all the examples, the random imperfection concerns the Young's moduli of the members. We consider the  cost function~\eqref{eq:optimisation_subeq_cost}, which is the weighted sum of the expectation and standard deviation of the structural compliance. We choose the respective weights in dependence on a single prescribed cost function parameter~$\alpha$ as 
\begin{equation}
\alpha_E = \frac{\alpha}{\overline J^*} \, , \quad  
\alpha_\sigma = \frac{1-\alpha}{\sigma^*_J} \,  . 
\end{equation}
The parameter~$\alpha$ is chosen such that~\mbox{$ 0 \le \alpha \le 1$.} The normalisation constant~$\overline J^*$ is the cost function value obtained by optimising with~$\alpha=1$, and the normalisation constant~$\sigma_J^*$ is the cost function value obtained with~$\alpha=0$. Note that the robust optimisation problem can also be interpreted as a multicriteria optimisation so that the obtained solution for a given~$\alpha$ represents a Pareto point~\citep{beyer2007robust}. We solve the discrete optimisation problem in all examples using the Method of Moving Asymptotes (MMA)~\citep{svanberg1987method}. 
%
\subsection{Verification example \label{sec:validation}}
\begin{figure}
	\centering
	\includegraphics[width=1\linewidth]{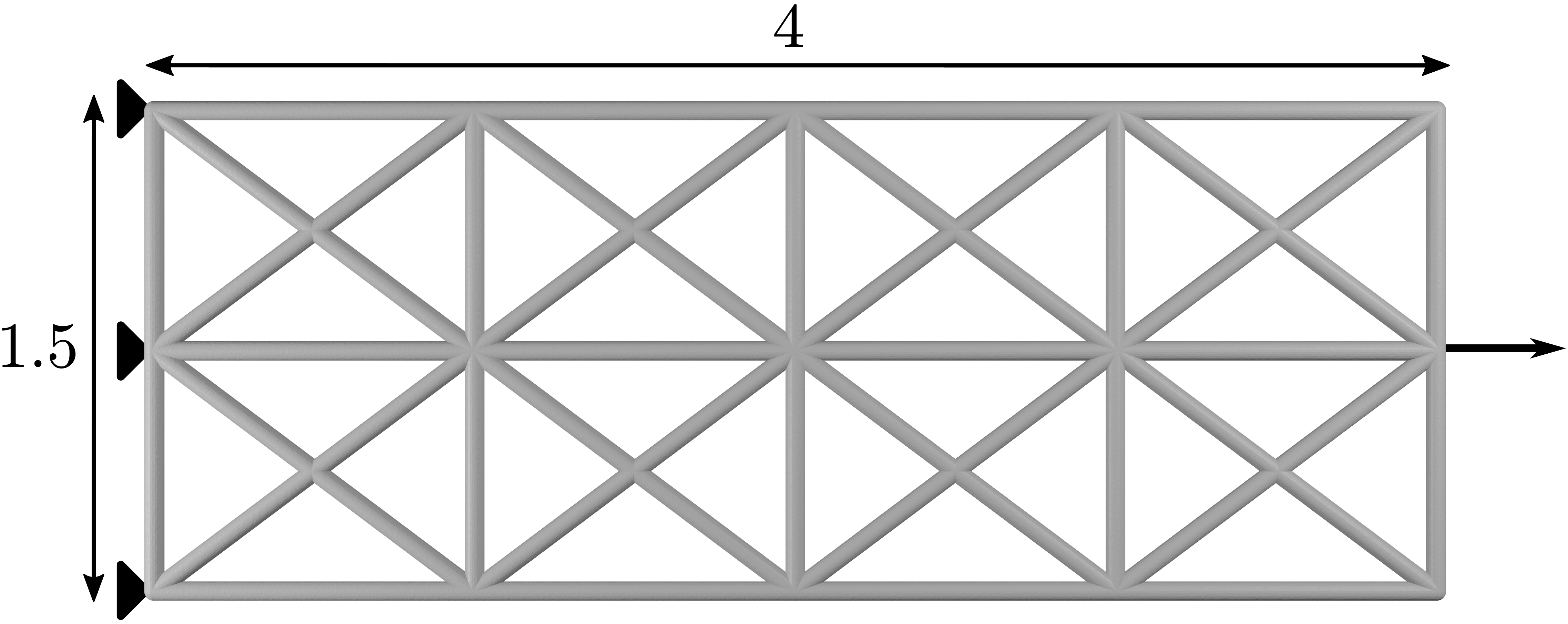}
	\caption{Verification example. Original structure, boundary conditions and loading.}
	\label{fig:asadground}
\end{figure}

%
To verify the correctness of our derivations and implementation, we consider the lattice structure depicted in Figure~\ref{fig:asadground}, which was previously investigated by~\cite{asadpoure2011robust}. The lattice consists of $2 \times 4$ unit cells of size $0.75 \times 1$. There are in total 15 joints and 38 members. The joints at the left end are fixed, and a single external force of magnitude~$1$ is applied at the centre of the right end.  The prescribed maximum volume is~\mbox{$V_{\max}=0.5$} and the cross-sectional member areas are constrained not to exceed~\mbox{$A_{\max} = 1$}. The Young's moduli of the members are random, spatially uncorrelated and have the probability density  
\begin{equation}
	p(\vec r) = \set N ( 100  \vec I_{n_e}, 10^2 \vec I_{n_e \times n_e} ) \,  , 
\end{equation}
where~$\vec I_{n_e}$ is an all-ones vector and~$ \vec I_{n_e \times n_e}$ an identity matrix. We note that the corresponding finite element solution has a non-diagonal covariance matrix as given in~\eqref{eq:solCovariance}.

\begin{figure}
	\centering
	\subfloat[$\alpha = 1$, $\overline{J} = 0.32$, $\sigma_J = 1.60\cdot 10^{-2}$.]{\includegraphics[width=0.85\linewidth]{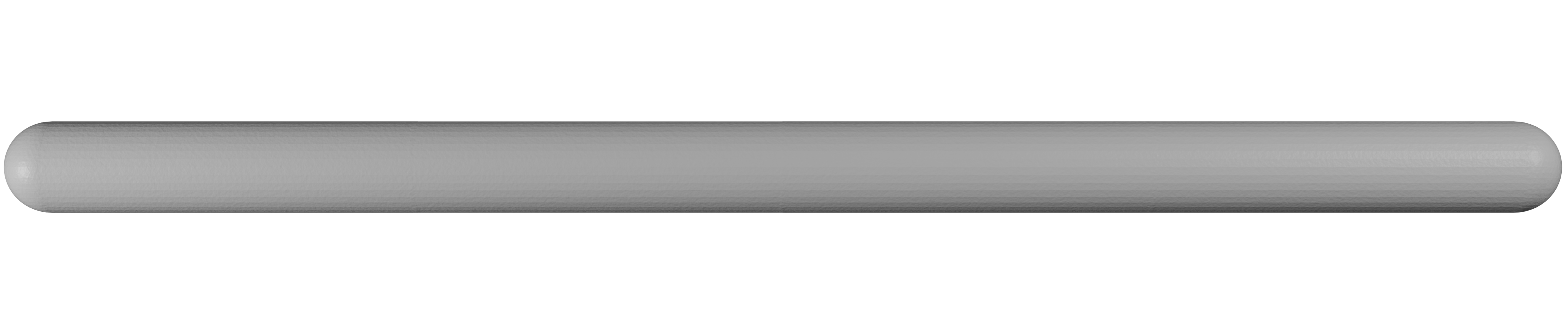}\label{fig:asad1}}\\
	\subfloat[$\alpha = 0.3$, $\overline{J} = 0.47$, $\sigma_J = 1.18\cdot 10^{-2}$.]{\includegraphics[width=0.85\linewidth]{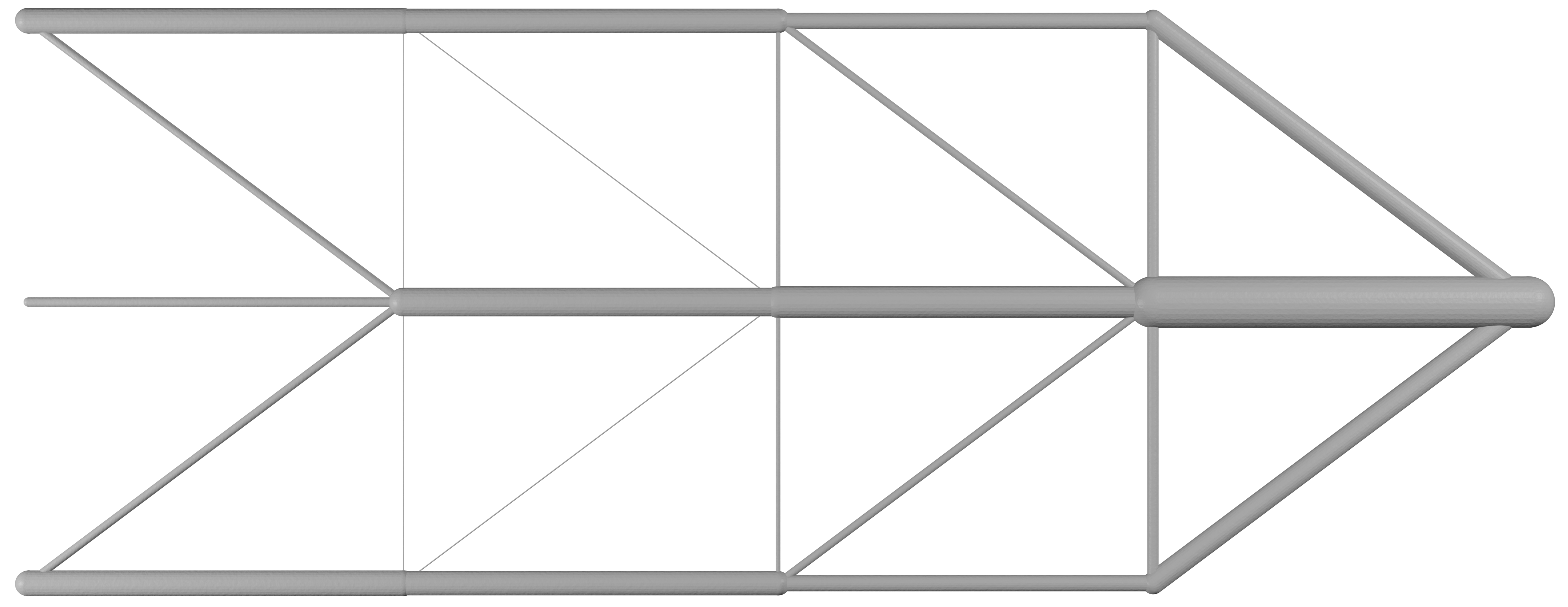}\label{fig:asad03}}\\
	\subfloat[$\alpha = 0$, $\overline{J} = 0.50$, $\sigma_J = 1.15\cdot 10^{-2}$.]{\includegraphics[width=0.85\linewidth]{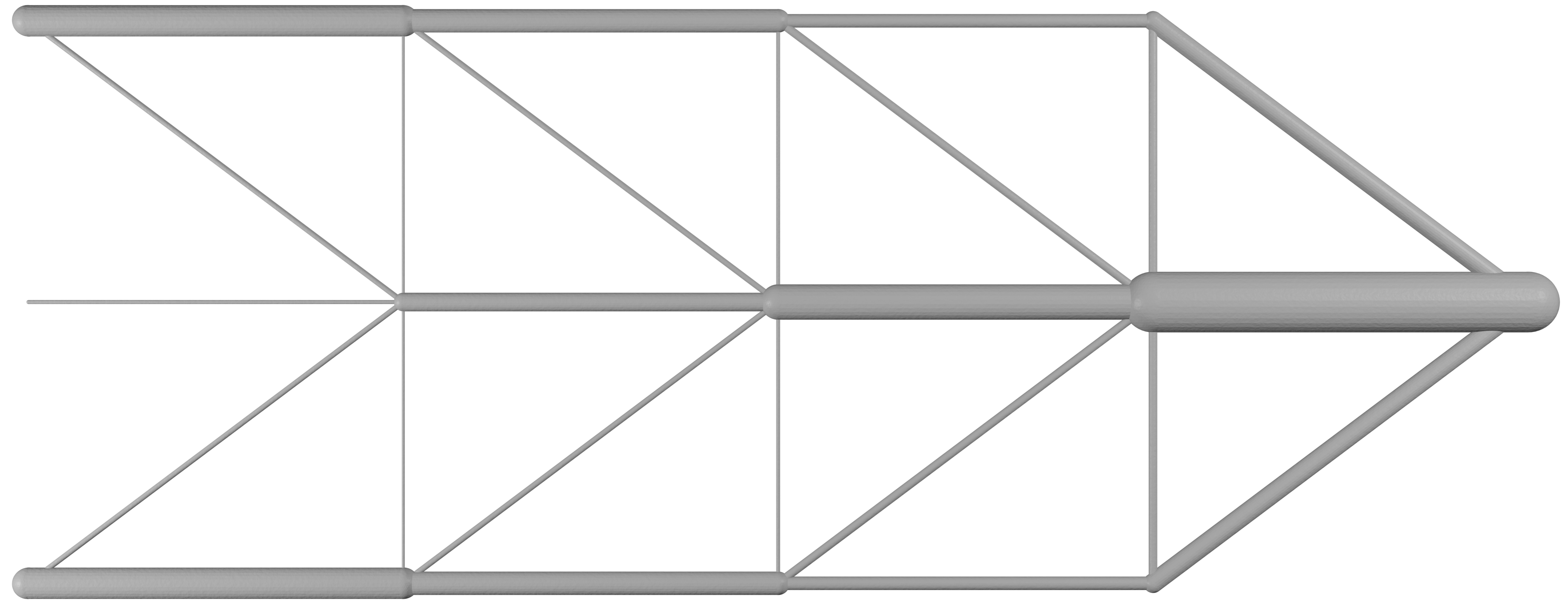}\label{fig:asad0}}
	\caption{Verification example. Optimised lattice structures for the cost function parameter values~$\alpha \in \{1, \, 0.3, \, 0 \}$. Optimising only with respect to the expectation of the compliance, i.e.~$\alpha=1$, yields an optimised structure consisting of only a colinear chain of four members. Including the standard deviation of the compliance by choosing $\alpha<1$ the members become more dispersed. For~$\alpha=0$ the cost function consists only of the standard deviation of the compliance. \label{fig:asadpoure}}
\end{figure}

When optimising only with respect to the expected compliance, i.e. by choosing $\alpha=1$, in the optimised structure shown in Figure~\ref{fig:asad1} only the members along the centre of the structure are present while all the others are missing. Note that this case is equivalent to conventional deterministic optimisation. Evidently, the obtained structure cannot be considered as robust given that there is only a single load path, and the compliance depends strongly on the random Young's moduli of the four members. 

To improve robustness, we take into account the standard deviation of the compliance by choosing a smaller~$\alpha$ value. As expected, the optimised structure  for~$\alpha=0.3$ depicted in Figure~\ref{fig:asad03} has many more members, a slightly higher compliance and a significantly smaller standard deviation as desired. Irrespective of its practicality, it is feasible only to minimise the standard deviation of the compliance while neglecting its expected value by choosing~$\alpha=0$. The respective optimised structure is shown in Figure~\ref{fig:asad0}. It bears emphasis that the material volume is the same in all the three optimised structures in Figure~\ref{fig:asadpoure}. Finally, we confirm that the obtained optimised structures are in close agreement with the results reported in~\cite{asadpoure2011robust}, verifying the correctness of our implementation. 
%
\subsection{Tension strip \label{sec:pillar}}
\begin{figure}
	\centering
	\includegraphics[width=\linewidth]{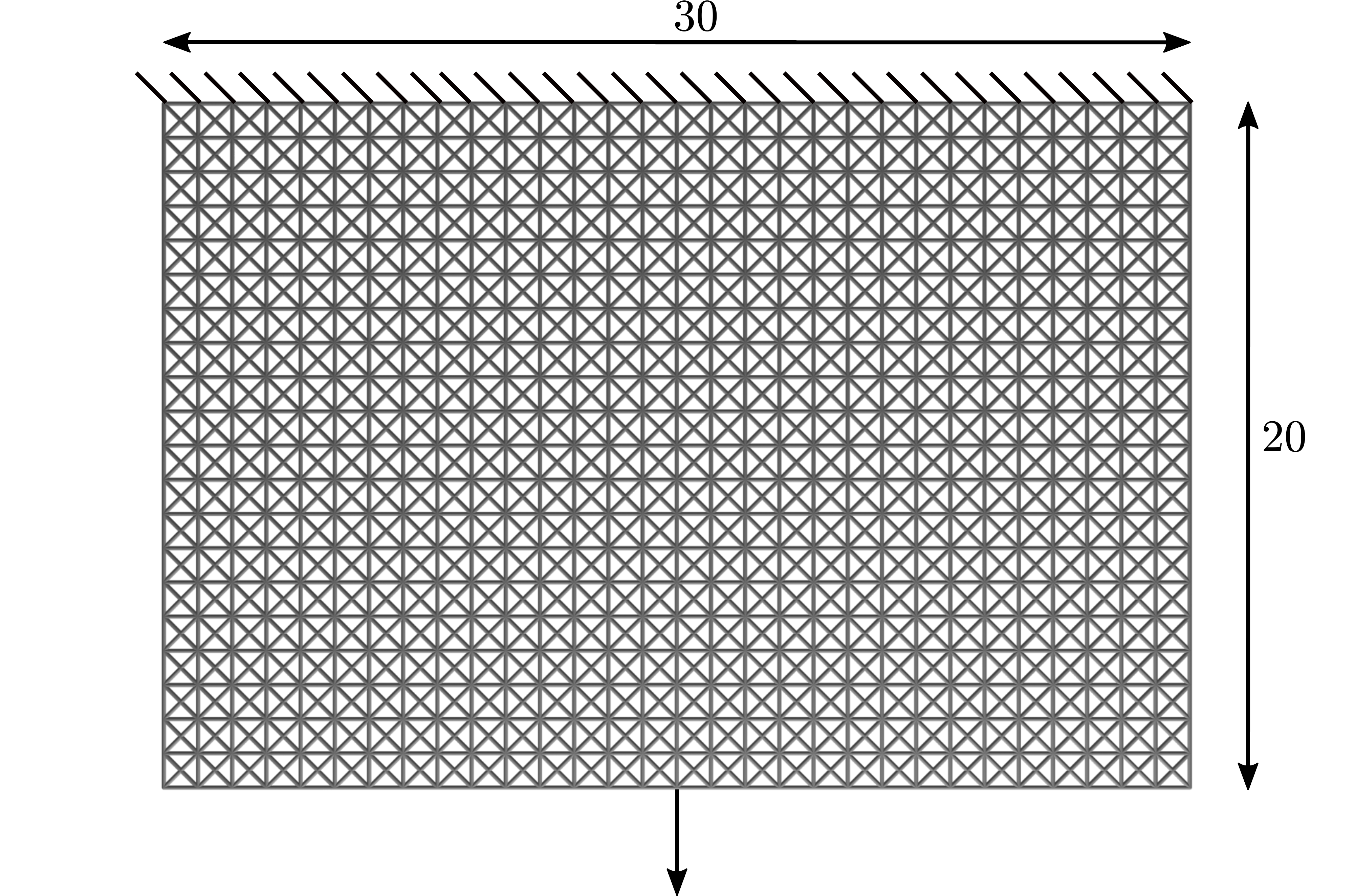}
	\caption{Tension strip. Original structure, boundary conditions and loading \label{fig:pillarground}}
\end{figure}

%
In this example, we optimise the two-dimensional tension strip depicted in Figure~\ref{fig:pillarground} and study the effect of the spatial properties of the random Young's moduli. The precision matrix of the respective probability density function is  determined using the SPDE formulation of random fields and its extension to lattice structures introduced in Section~\ref{sec:maternLattices}.  The width and height of the tension strip are $30 \times 20$. The joints along the top boundary are fixed, and all the remaining boundary joints are unconstrained. A single load of magnitude~$1$ is applied at the centre on the bottom edge. Each unit cell is of size~$1 \times 1$ and has two diagonal members. The structure has 2450 members, and the mean of Young's modulus for all members is~$\overline E = 100$. Furthermore, the total allowable volume is $V_{\max} = 6$ and the cross-sectional member areas are constrained not to exceed $A_{\max} = 1$. The radius of the density filter is chosen as $R = 2.5$, i.e. $2.5$  times the unit cell size. 
%
\subsubsection{Isotropic, stationary covariance \label{subsec: isohomo}}
\begin{figure}
	\centering
	\subfloat[$\ell =0.5$]{\includegraphics[width=0.31\textwidth]{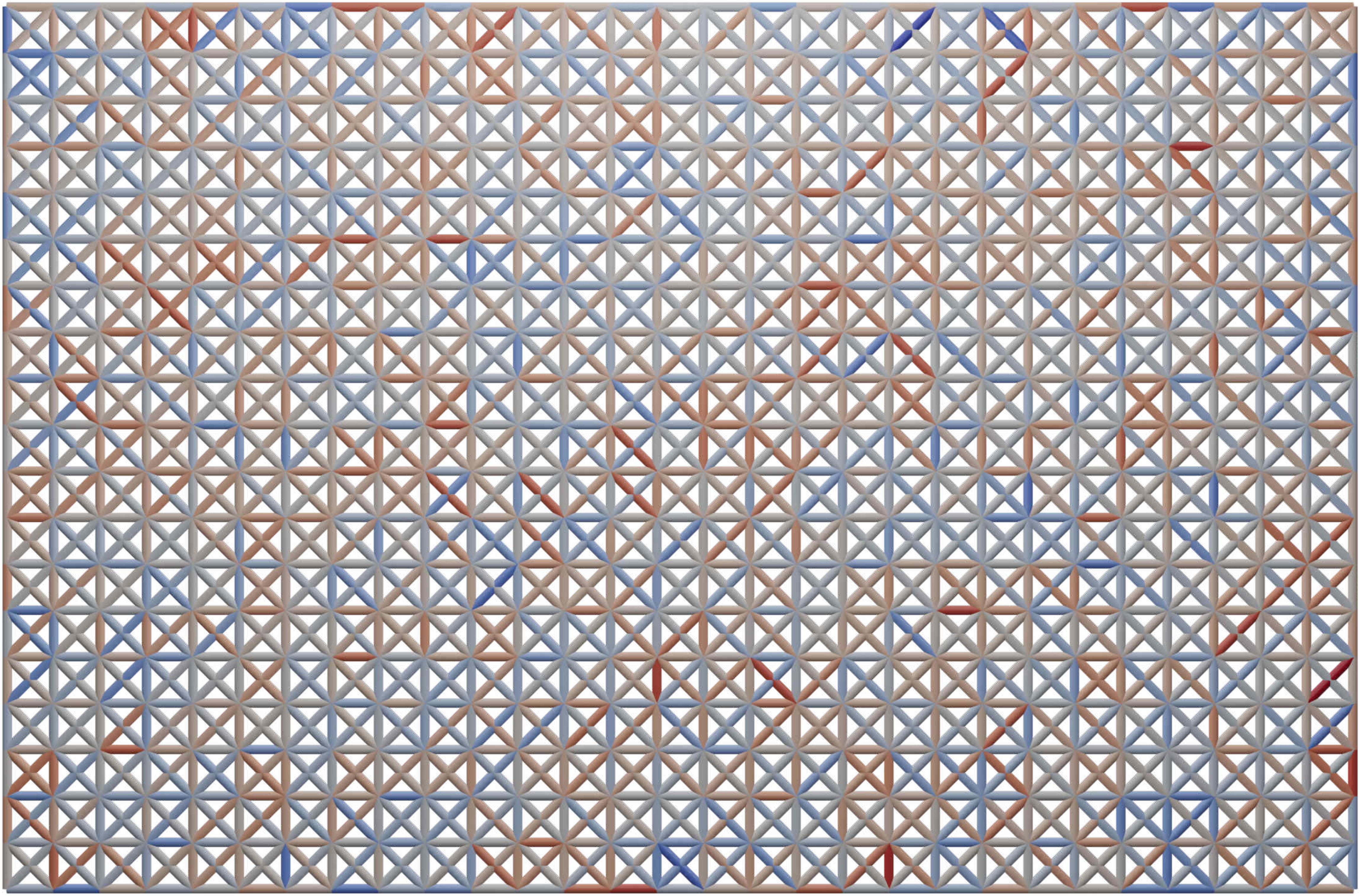}\label{fig:pillar_sample_05}}
	\\
	\subfloat[$\ell =5$]{\includegraphics[width=0.31\textwidth]{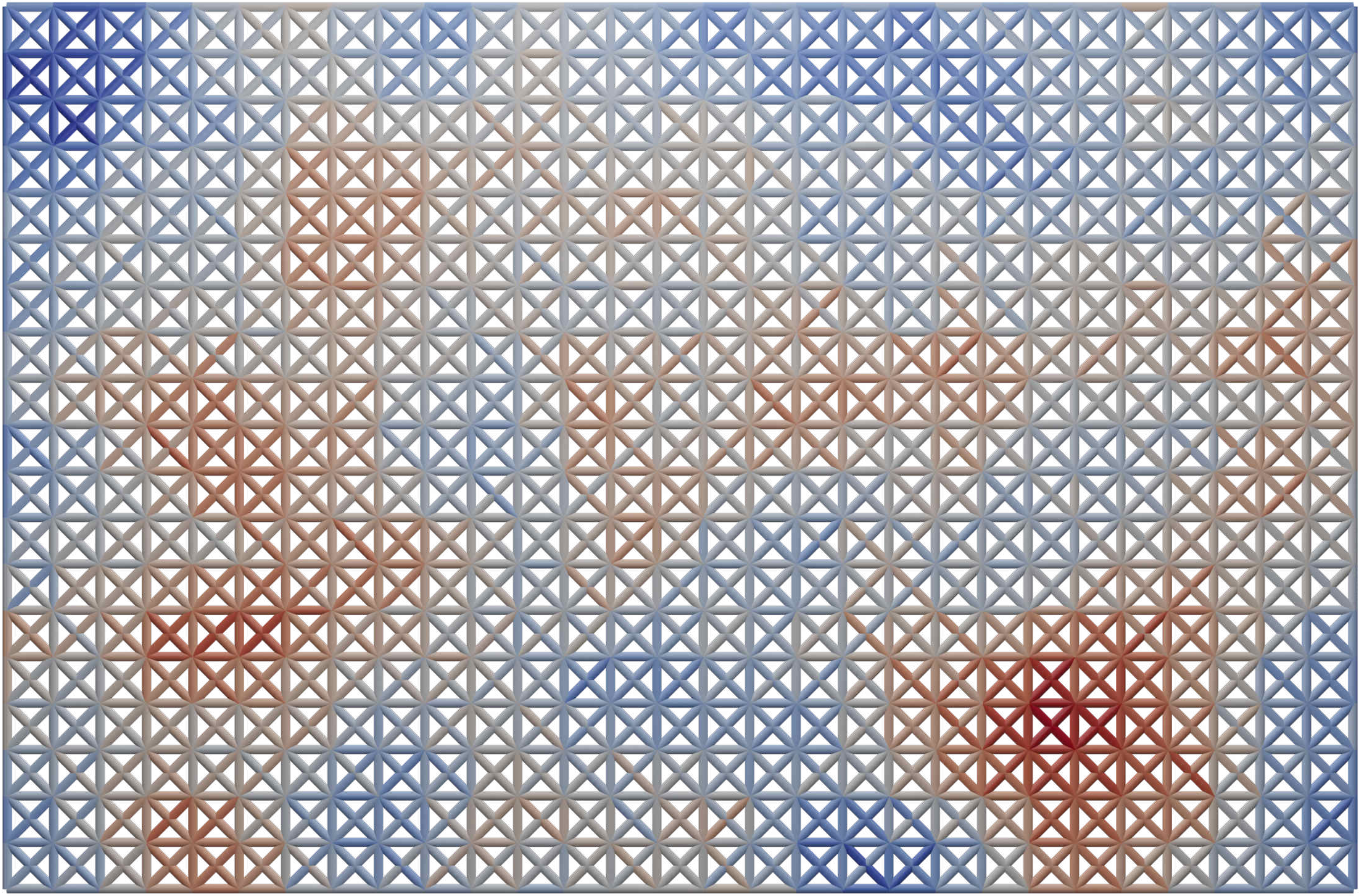}\label{fig:pillar_sample_5}}
	\\
	\subfloat[$\ell =10$]{\includegraphics[width=0.31\textwidth]{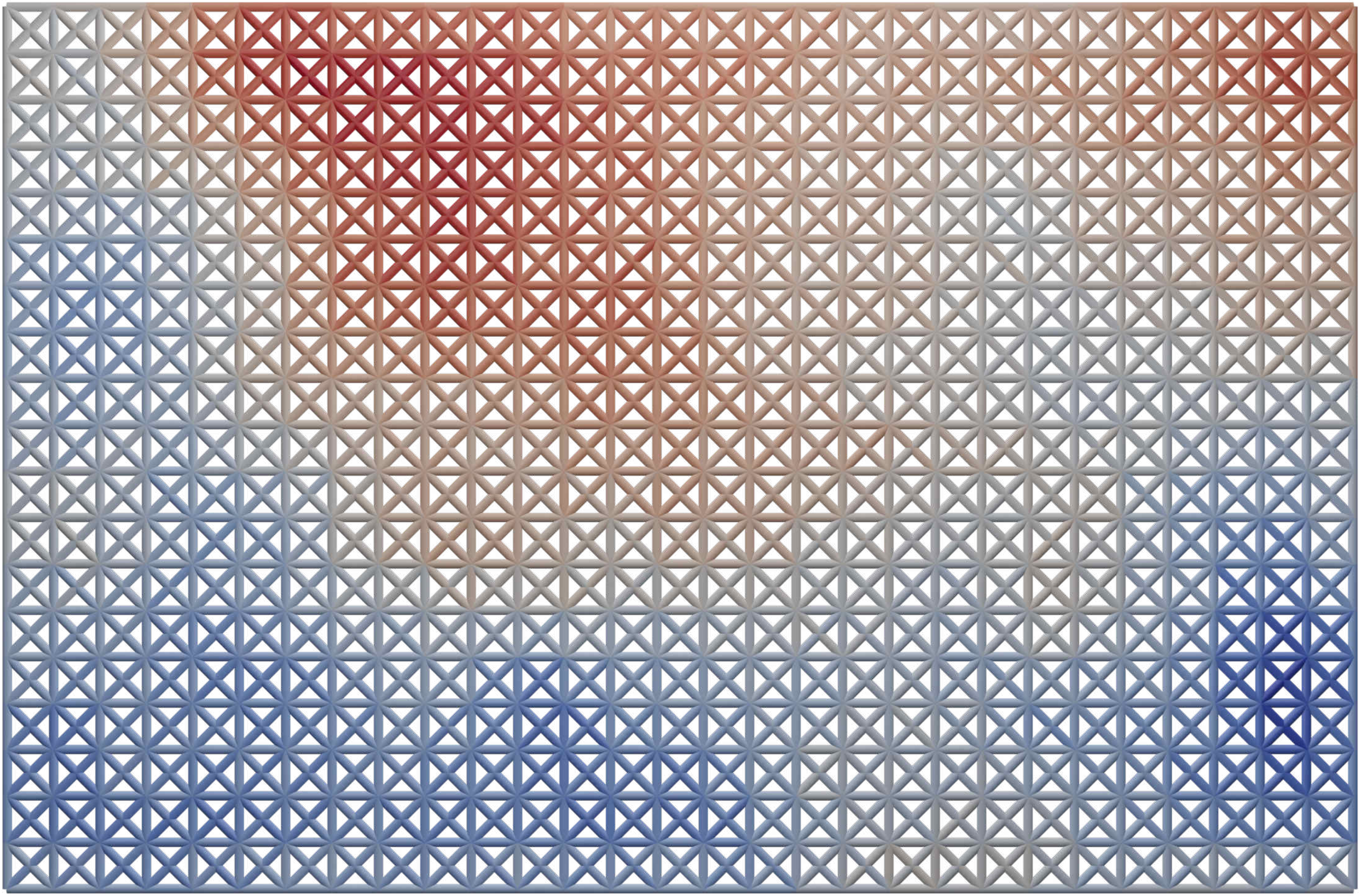}\label{fig:pillar_sample_10}}
	\caption{Tension strip with isotropic, stationary covariance. Isocontours of the member Young's moduli obtained by sampling from the SPDE-defined Gaussian probability density with a covariance length-scale~\mbox{$\ell \in \{ 0.5, 5, 10 \}$}, smoothness~$\nu =1.5$, standard deviation~ $\sigma = 10$ and constant mean $\overline E = 100$.}
    \label{fig:pillar_iso_samples}
\end{figure}

%
First, we investigate the effect of the covariance length-scale \mbox{$\ell \in \{ 0.5, 5.0, 10.0 \}$} and the cost function weight parameter \mbox{$\alpha \in \{ 1, \, 0.5, \, 0 \}$} on the optimisation results. The covariance smoothness parameter and the standard deviation are \mbox{$\nu = 1.5$} and \mbox{$\sigma = 10$}.  For illustration purposes, in Figure~\ref{fig:pillar_iso_samples}, three samples from the considered random fields with different length-scales are shown. Note that~\mbox{$E(\vec x) - \overline E$} is equal to the solution of the SPDE~\eqref{eq:spdeFErecurs}  with a Gaussian white noise forcing. We use this fact to sample efficiently large random fields like the ones depicted. In Figure~\ref{fig:pillar_iso_samples},  the members are coloured according to their Young's modulus, equal to the solution of the SPDE plus the mean~$\overline E$. For~$\ell \lesssim 1$, the Young's moduli of the members are effectively uncorrelated for the considered lattice with a unit cell size~$1 \times 1$.

The nine optimised lattice structures for different combinations of length-scale parameters~$\ell$ and cost function weight parameters~$\alpha$ are depicted in Figure~\ref{fig:homo_iso}. When the cost function is only the expectation of the compliance, i.e.~$\alpha=1$, the optimal structure and the respective expected compliance~$\overline J$ do not depend on the length-scale. However, an increase in~$\ell$ still leads to a rise in the standard deviation of the compliance~$\sigma_J$. The cost function depends on the standard deviation when choosing~$\alpha<1$. As seen in the middle and right columns in Figure~\ref{fig:homo_iso}, this leads to designs with a more spread-out member layout.  The extent of the spreading is inversely proportional to the length-scale~$\ell$. It is  intuitively clear that a decrease in the length-scale leads  to increasingly uncorrelated Young's moduli, in which case the  best robust optimisation strategy consists of distributing the load to as many as possible members while satisfying all the constraints.  
\begin{figure*}
	\centering
	\includegraphics[width=\linewidth]{homo_iso.pdf}
	\caption{Tension strip with isotropic, stationary covariance. Optimised structures for different combinations of length-scale~$\ell$ and cost function parameter~$\alpha$.  For each structure the minimum expected compliance~$\overline{J}$  and standard deviation of the compliance~$\sigma_J$  are given at the bottom left of each subfigure. \label{fig:homo_iso}}
\end{figure*}

%
\subsubsection{Non-stationary, anisotropic covariance \label{subsec: nonstandard}}
\begin{figure}[t]
	\centering
	\subfloat[]{\includegraphics[scale=0.09]{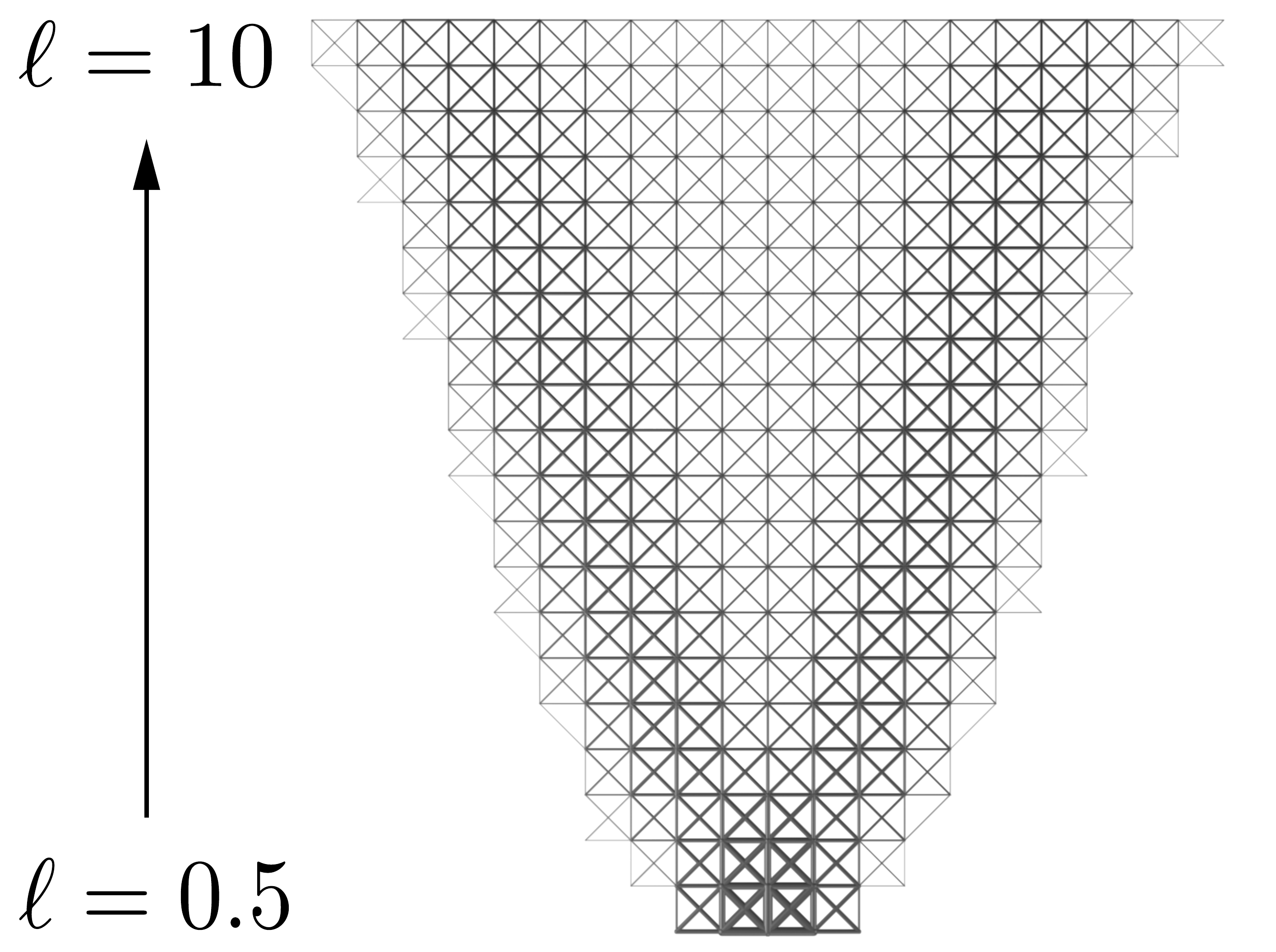}\label{fig:nonhomo_a}}  
	\\
	\subfloat[]{\includegraphics[scale=0.09]{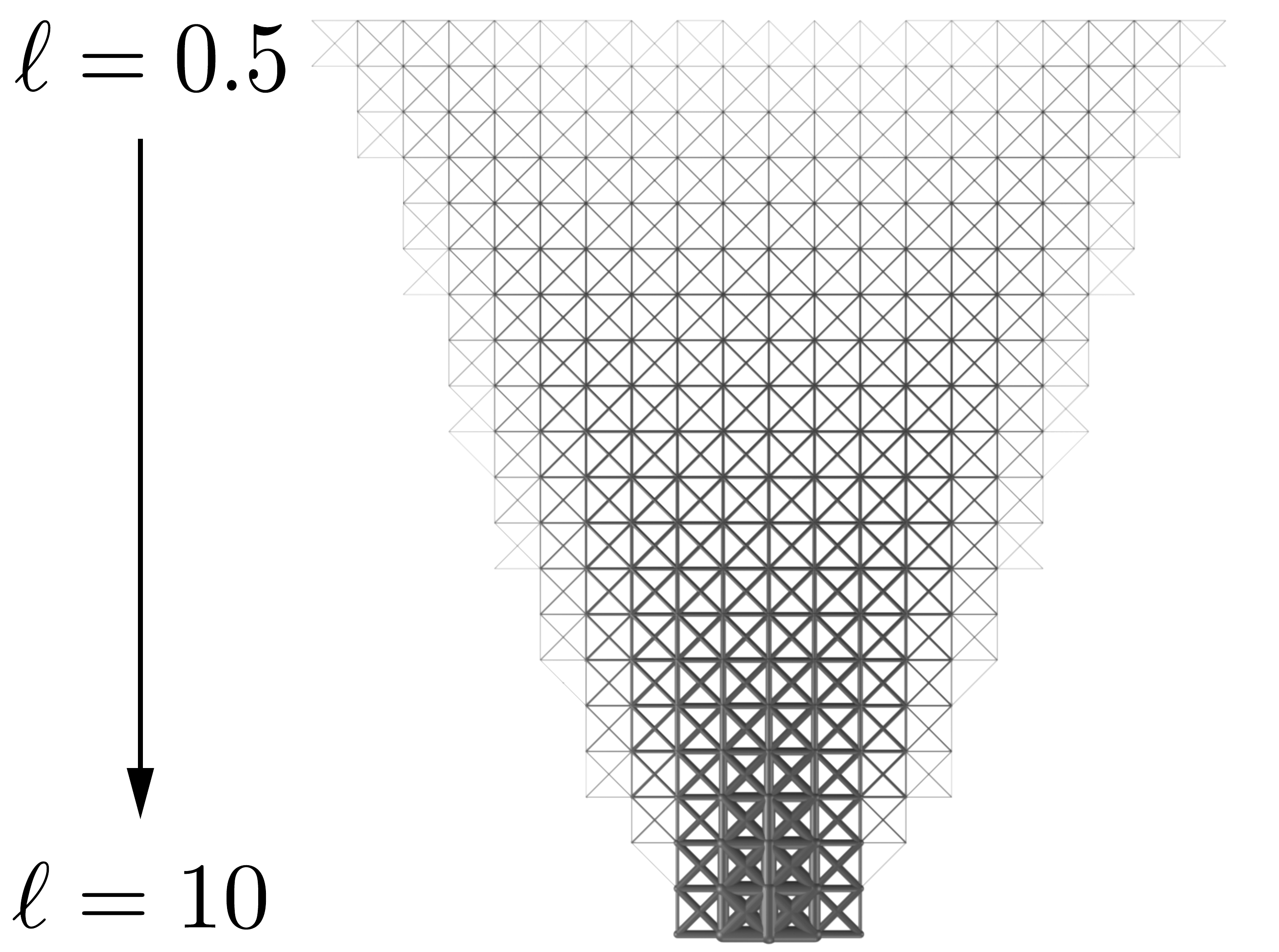}\label{fig:nonhomo_b}}  
	\caption{Tension strip with non-stationary, isotropic covariance.  In (a) the length-scale increases linearly from 0.5 (bottom) to 10 (top). In contrast, in (b) the length-scale increases linearly from 0.5 (top) to 10 (bottom).\label{fig:nonhomo}}
\end{figure}

%
We now optimise a tension strip with a random Young's modulus with a non-stationary isotropic covariance with a length-scale $\ell(\vec x)$ varying linearly from the top to the bottom. The obtained optimised structures for two different length-scale distributions are compared in Figure~\ref{fig:nonhomo}.
In Figure~\ref{fig:nonhomo_a}, the length-scale increases linearly from the bottom to the top according to  \mbox{${\ell(y) = 0.5(0.95 y+1)}$}, so that the length-scale at the top is \mbox{${\ell(20)=10}$} and at the bottom ${\ell(0)=0.5}$. In contrast, in Figure~\ref{fig:nonhomo_b}, the length-scale increases linearly from the top to the bottom according to ${\ell(y) = 0.5(-0.95 y+20),}$ leading to ${\ell(20)=0.5}$ and ${\ell(0)=10}$. The other covariance parameters ${\sigma=10}$ and ${\nu = 1.5}$ are the same for both structures. The Young's moduli vary more rapidly in shorter length-scale regions, prompting the algorithm to spread the material among more members.

\begin{figure}
	\centering
	\subfloat[$d_{\vert\vert}= 1,\,d_{\perp}= 5$]{\includegraphics[scale=0.085]{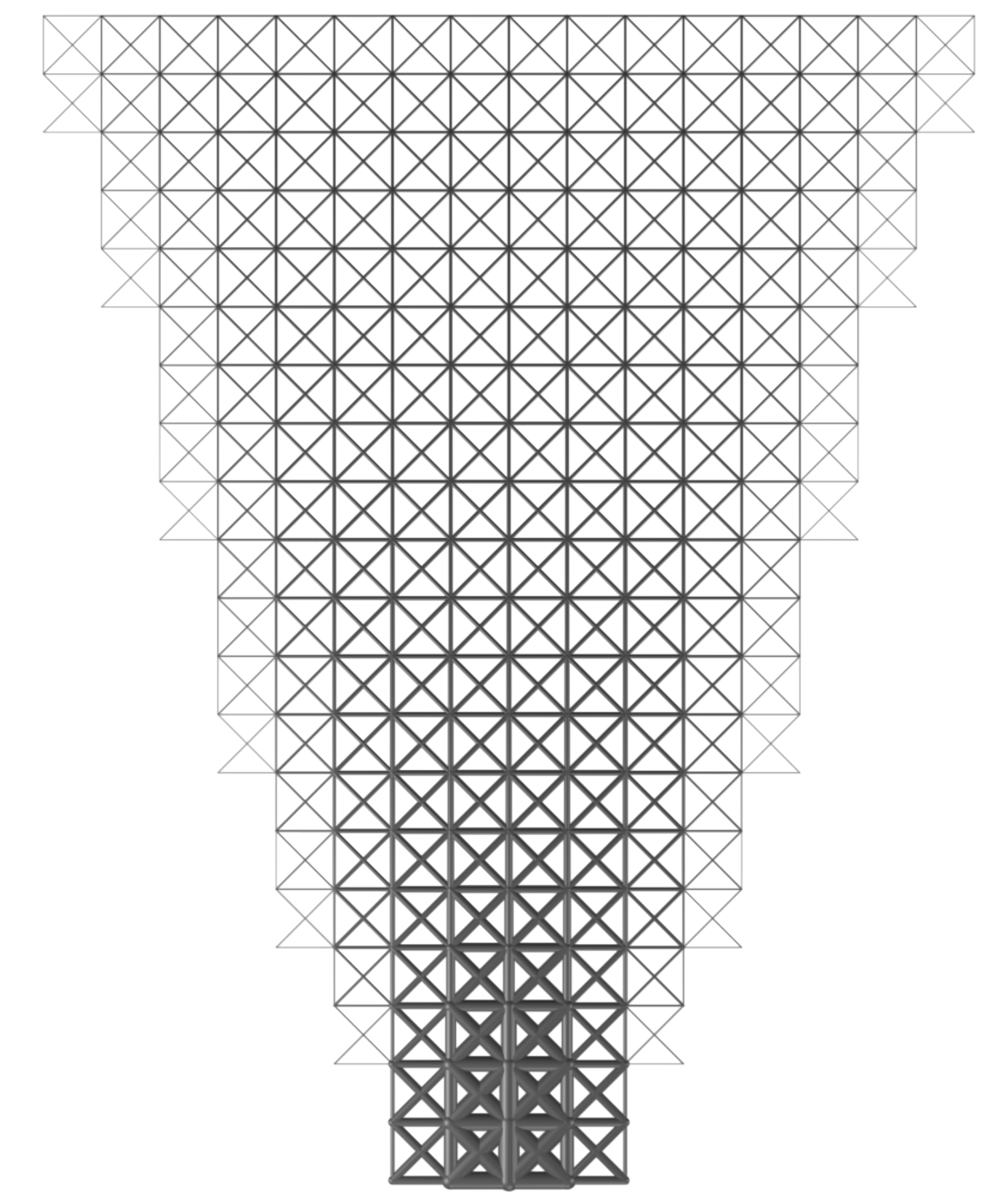}\label{fig:aniso_a}}  \hspace{0.05\textwidth}
	\subfloat[$d_{\vert\vert}= 5,\,d_{\perp}= 1$]{\includegraphics[scale=0.085]{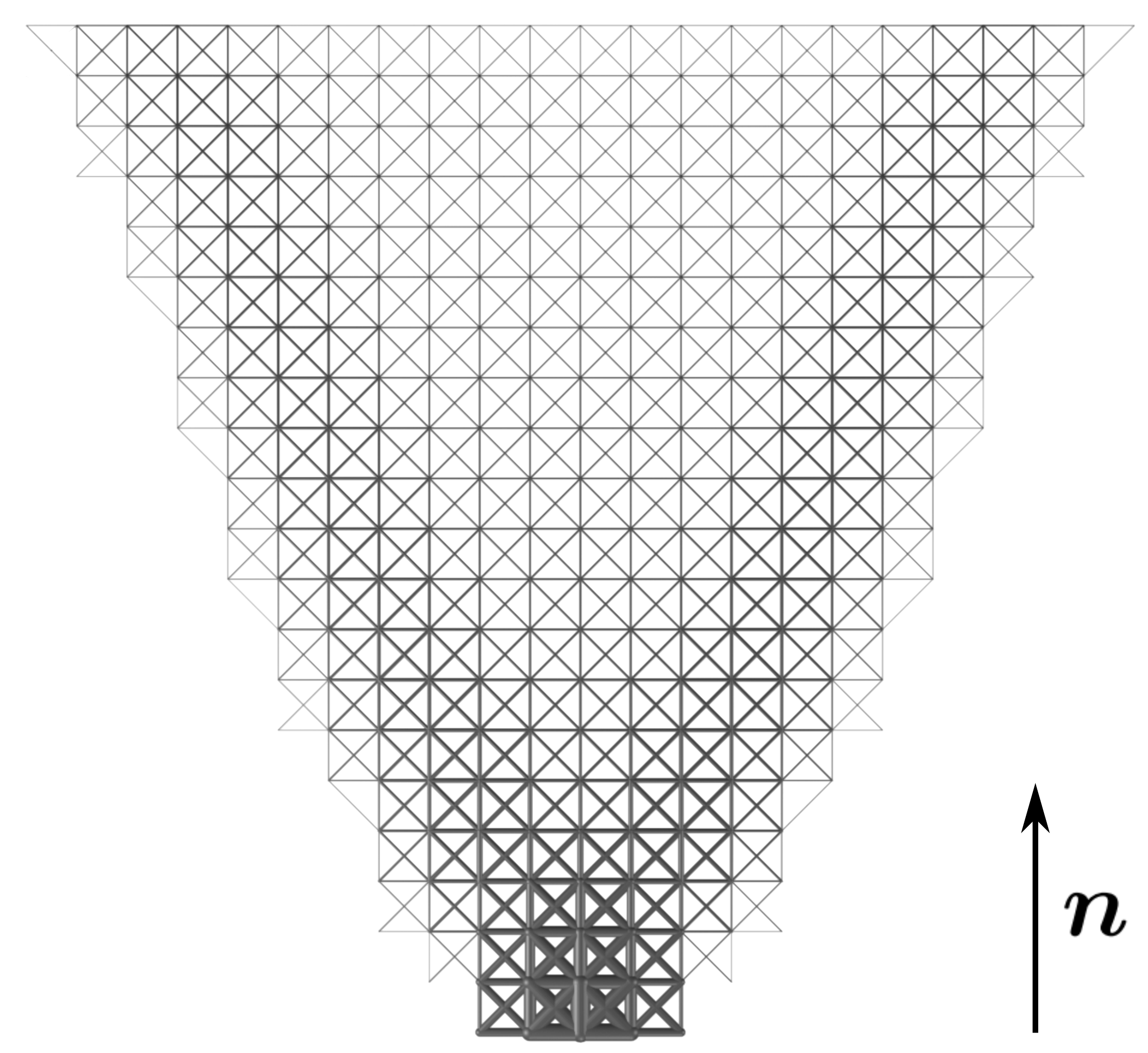}\label{fig:aniso_b}}
	\caption{Tension strip with stationary, anisotropic covariance.  Final optimised structures using stationary, anisotropic random Young's moduli. 
	 \label{fig:aniso} }
\end{figure}

At last, we consider a tension strip with two different anisotropic random Young's modulus distributions, see Figure~\ref{fig:aniso}. The prescribed unit vector~$\vec n$ gives the direction of the anisotropy. In the first case (Figure~\ref{fig:aniso_a}) with $d_{\vert\vert}= 1$ and $d_{\perp}= 5$   the members perpendicular to~$\vec{n}$ have higher uncertainty. The robust optimisation leads to a narrower structural system as the load is carried primarily by the vertical members with a lower uncertainty.  In contrast, in the second case (Figure~\ref{fig:aniso_b}) with $d_{\vert\vert}= 5$ and $d_{\perp}= 1$ the load is carried mainly by the diagonal members due to the high uncertainty associated with the vertical members. 

\newif\ifcantilever
\cantilevertrue

\ifcantilever

%
\subsection{Cantilever\label{subsec: penal}}
%
We consider next a cantilever lattice structure of size $20 \times 10$ with unit cells of size $0.5 \times 0.5$, see Figure~\ref{fig:penalground}. Each unit cell has two diagonal members, so there are 861 joints and 3260 members in total.  A vertical  load of magnitude 100 is applied at the mid-height on the right end of the lattice while the joints along the left end are fixed. The uncorrelated member Young's moduli have a constant mean value of \mbox{$\overline{E}=7\cdot 10^7$} and a standard deviation of~\mbox{$\sigma = 7\cdot 10^6$}. The prescribed maximum volume is $V_{\max}=1.6$, and the cross-sectional areas are constrained not to exceed \mbox{$A_{\max} = 5.1\times 10^{-3}$}. The density filter radius is chosen as \mbox{$R=1$}.
\begin{figure}
	\centering
	\includegraphics[width=\linewidth]{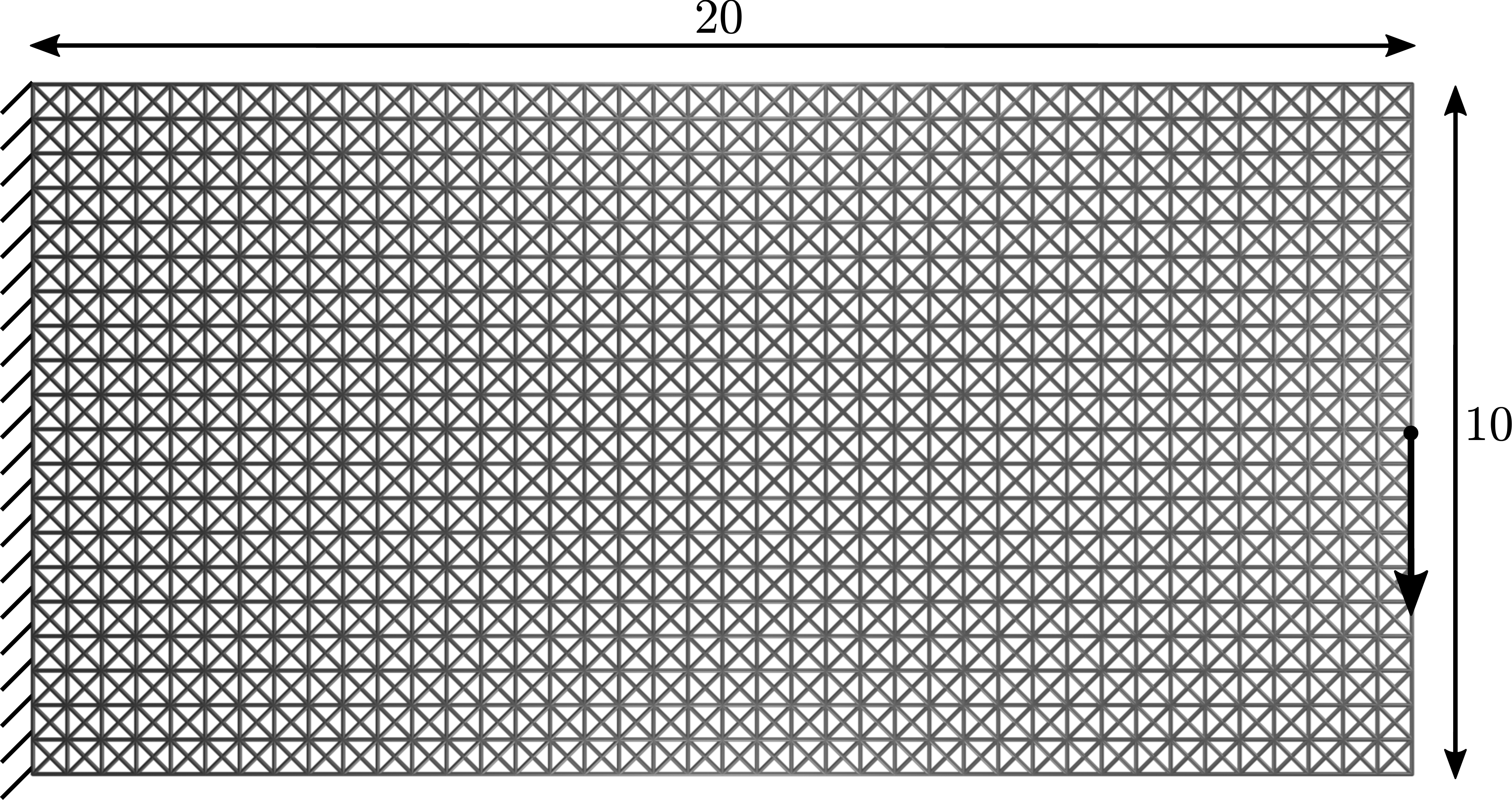}
	\caption{Cantilever. Original structure, boundary conditions, and loading.}
	\label{fig:penalground}
\end{figure}

In this example, we penalise the relative member densities $\tilde s_e (s_e)$ with the penalisation function introduced in Figure~\ref{fig:spline_plots} to obtain a manufacturable structure with cross-sectional areas larger than a prescribed threshold. The optimised structures in Figure~\ref{fig:penalisation} correspond to cost function weight parameter values \mbox{$\alpha \in \{ 1, \, 0.5, \, 0.1 \}$} and computations with and without penalisation. Selected histograms of the respective cross-sectional areas are given in Figure~\ref{fig:histogram_penal}. It is evident that in the case of no penalisation, there are a large number of, likely impossible to manufacture, members with very small cross-sectional areas. 

Furthermore, although the shape of the optimised structures is almost insensitive to~$\alpha$, the mean and standard deviation of the compliance are sensitive to the choice of~$\alpha$.  Figure~\ref{fig:pareto_curves} shows the respective Pareto fronts obtained by varying~$\alpha$. Without going into details, Figure~\ref{fig:pareto_curves}  also contains a third Pareto front for a penalisation function which is slightly less curved than the one depicted in Figure~\ref{fig:spline_plots}.     
\begin{figure*}
	\centering
	\includegraphics[width=\linewidth]{Penalisation.pdf}
	\caption{Cantilever. Comparison of optimised structures with and without penalisation of the relative member densities for different cost function weight parameter values.}
	\label{fig:penalisation}
\end{figure*}
\begin{figure*}
	\centering
	\subfloat[No penalisation, $\alpha = 1$]{\includegraphics[width=0.25\textwidth]{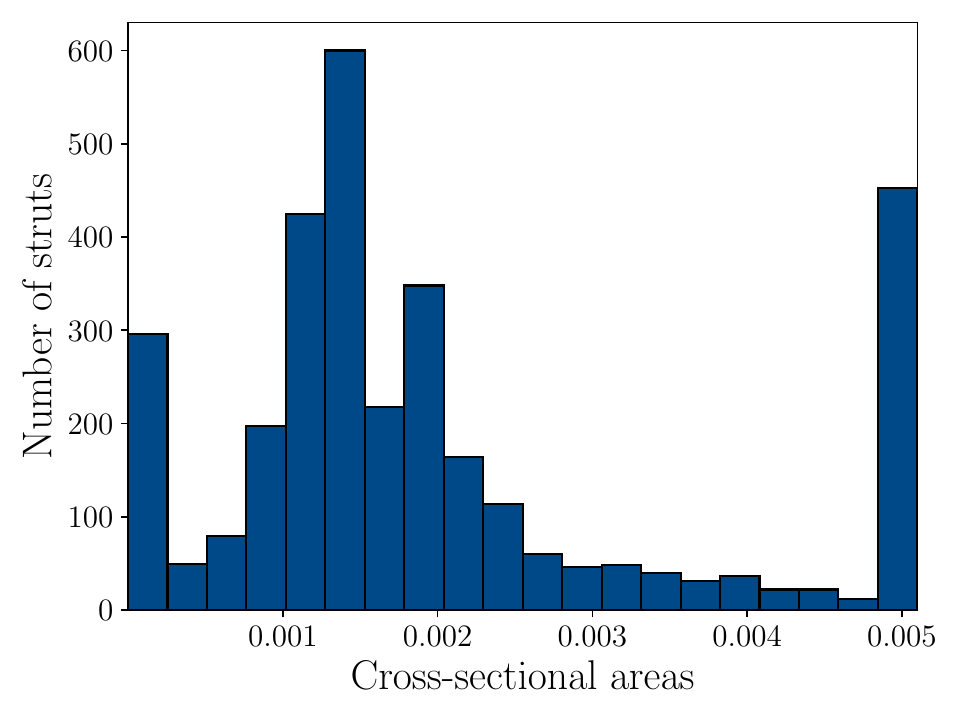}\label{histogram-p=1-a=1.00}}
	%
	%
        \subfloat[No penalisation, $\alpha = 0.50$]{\includegraphics[width=0.25\textwidth]{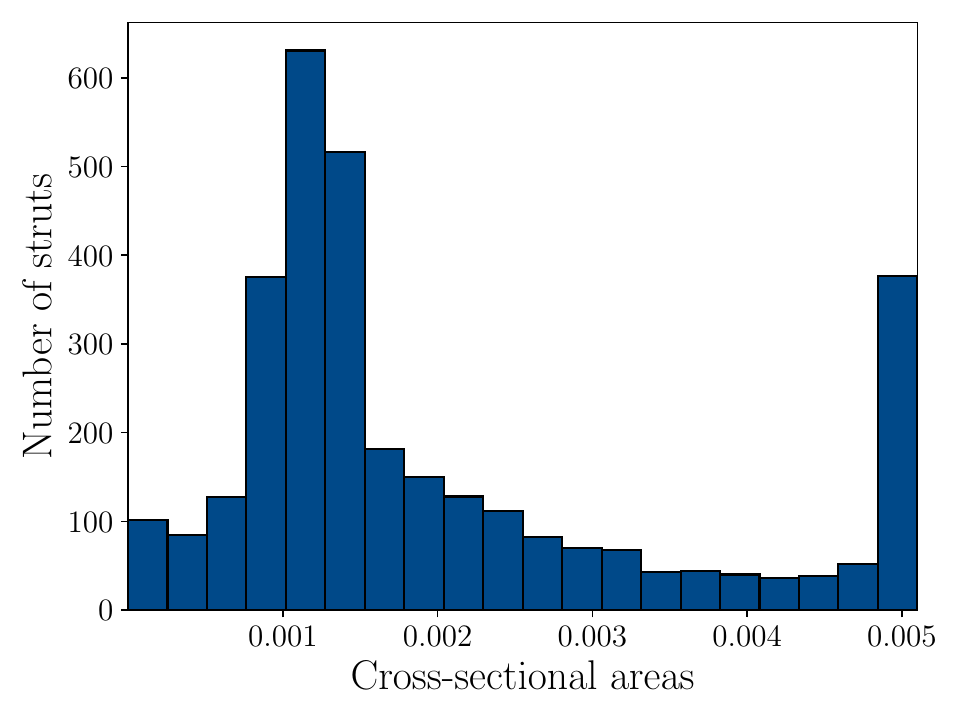}\label{histogram-p=1-a=0.50}}
	%
	\subfloat[With penalisation, $\alpha = 1$]{\includegraphics[width=0.25\textwidth]{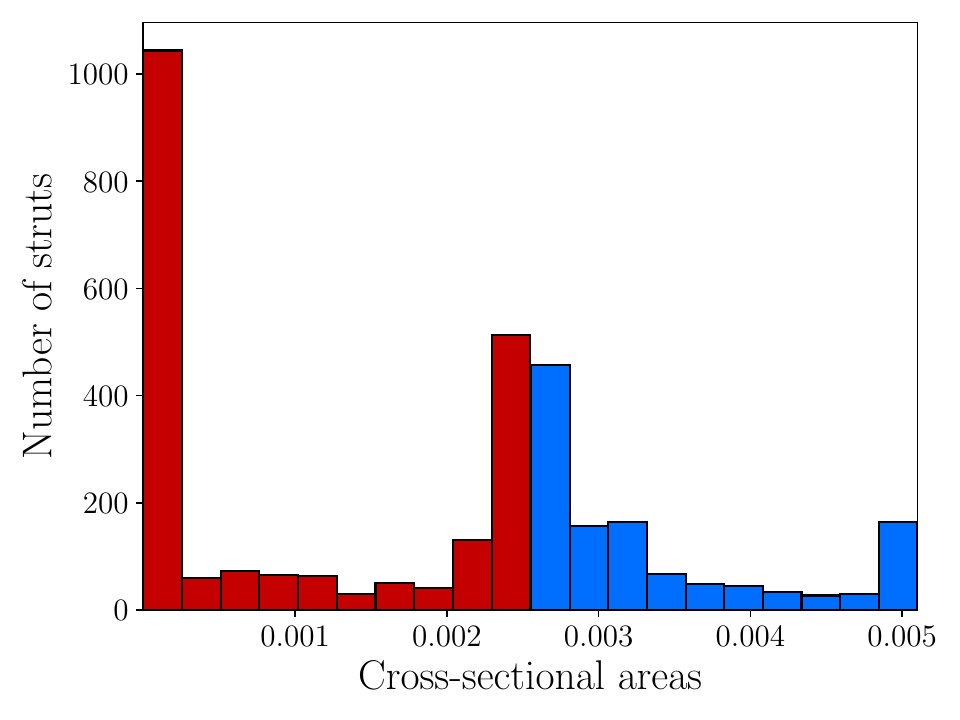}\label{histogram-B-spline-a=1.00}}
	%
	%
        \subfloat[With penalisation,  $\alpha = 0.50$]{\includegraphics[width=0.25\textwidth]{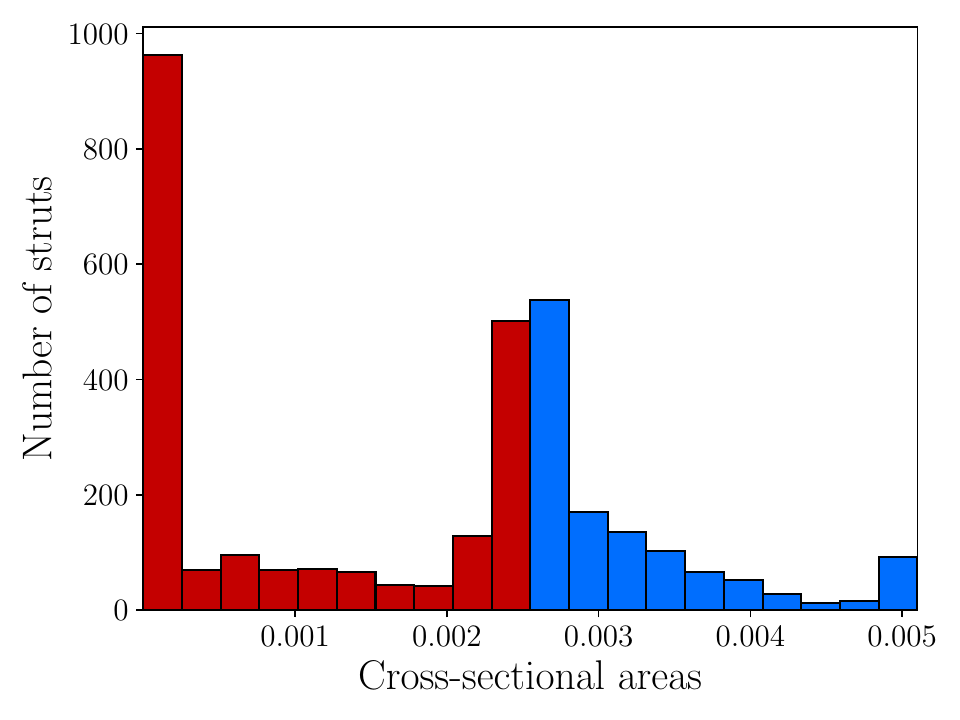}\label{histogram-B-spline-a=0.50}}
	%
	\caption{Cantilever. Selected histograms of the cross-sectional areas of the optimised structures shown in Figure~\ref{fig:penalisation}. The first two, (a) and (b), are for optimisation without and the second two, (c) and (d), with penalisation of the relative member densities. In (c) and (d) the red bars represent the hard to manufacture members with a relative density~$s^*_e < 0.5$. }
	\label{fig:histogram_penal}
\end{figure*}
\begin{figure}
    \centering
    \includegraphics[width=\linewidth]{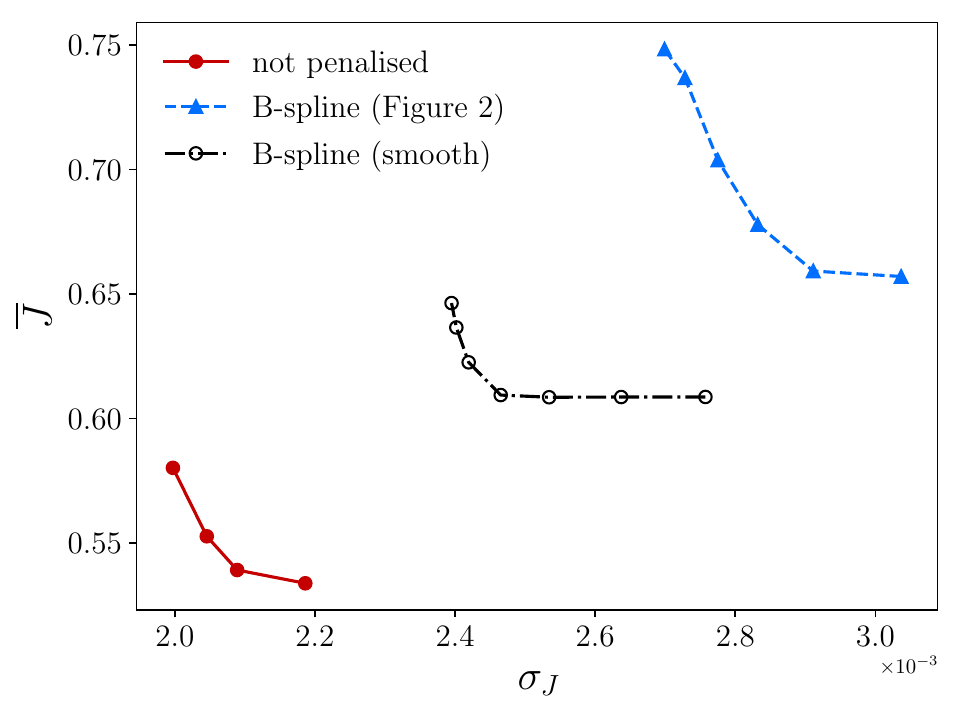}
    \caption{Cantilever. Pareto fronts of the multi-objective robust optimisation problem obtained by varying the cost function weight parameter~$\alpha$. The two B-spline penalised curves correspond to different penalisations of the relative member densities, i.e. according to Figure~\ref{fig:spline_plots} and a (not shown) less curved function.}
    \label{fig:pareto_curves}
\end{figure}

\fi

%
\subsection{Engine bracket \label{sec:bracket}}
\begin{figure}
	\centering
	%
	%
	\includegraphics[width=0.4\textwidth]{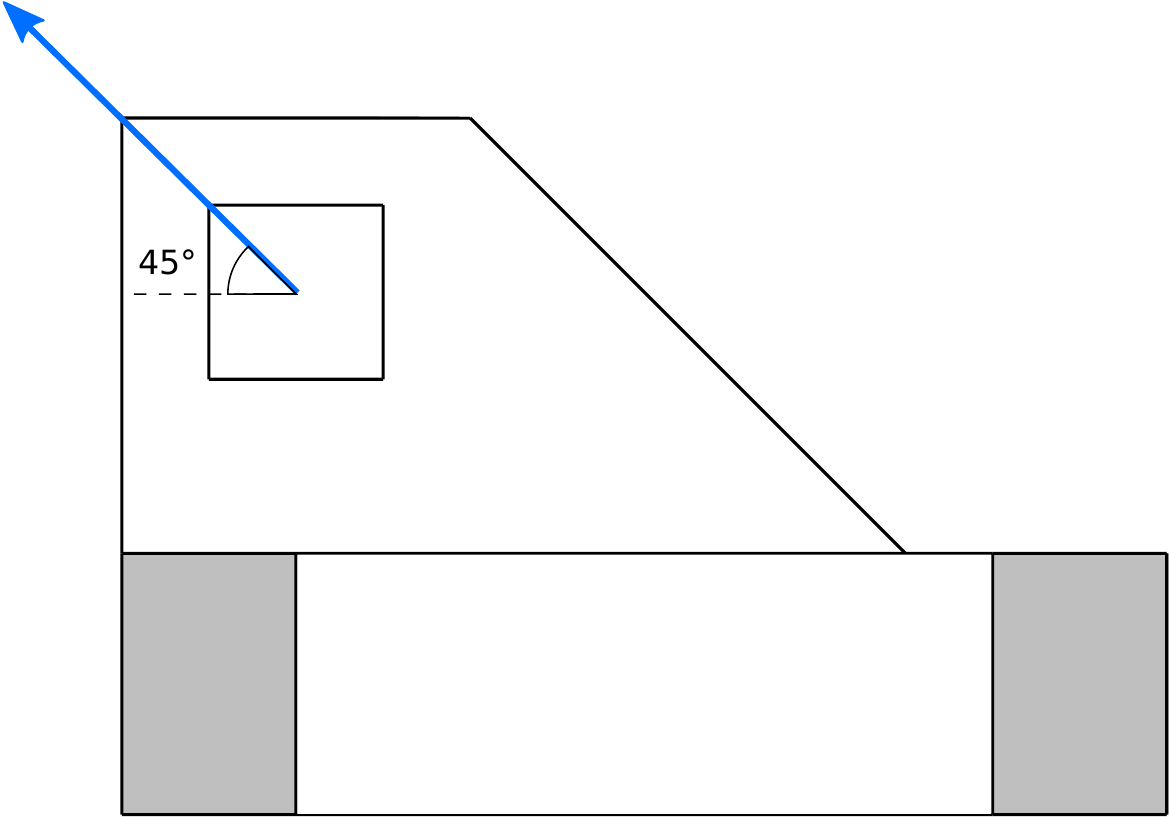}
	\caption{Engine bracket. Boundary conditions and loading. The arrow represents the load applied inclined by $45^\circ$ with respect to the base plate. All the  joints on the reentrant surfaces of the base plate (indicated in grey) are fixed.}
	\label{fig:bracketgroundside}
\end{figure}

\begin{figure}
	\centering
	\subfloat[Top view]{\includegraphics[width=0.5\textwidth]{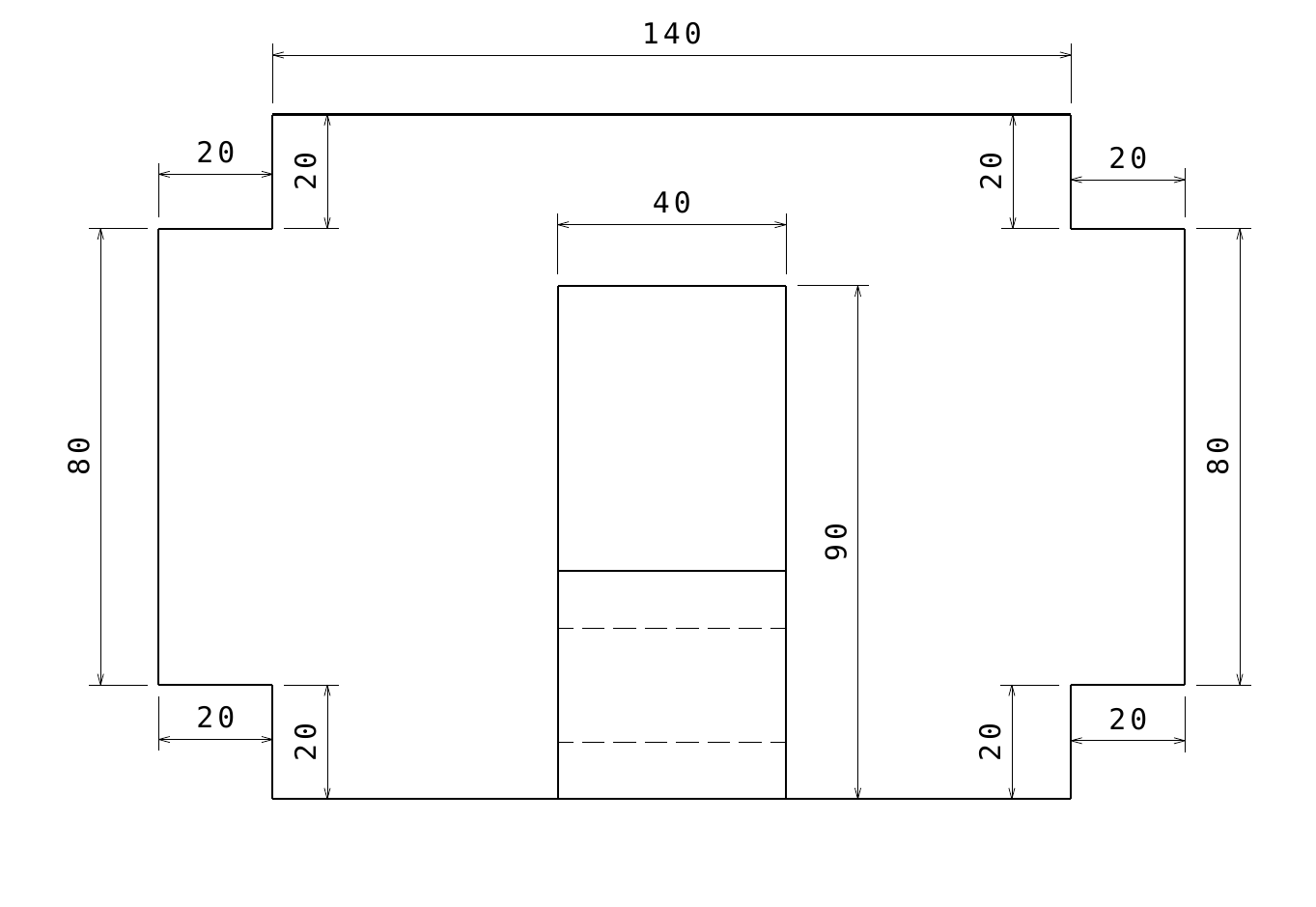}\label{fig:GE_bracket_tda}}
	\\
	\subfloat[Side view]{\includegraphics[width=0.5\textwidth]{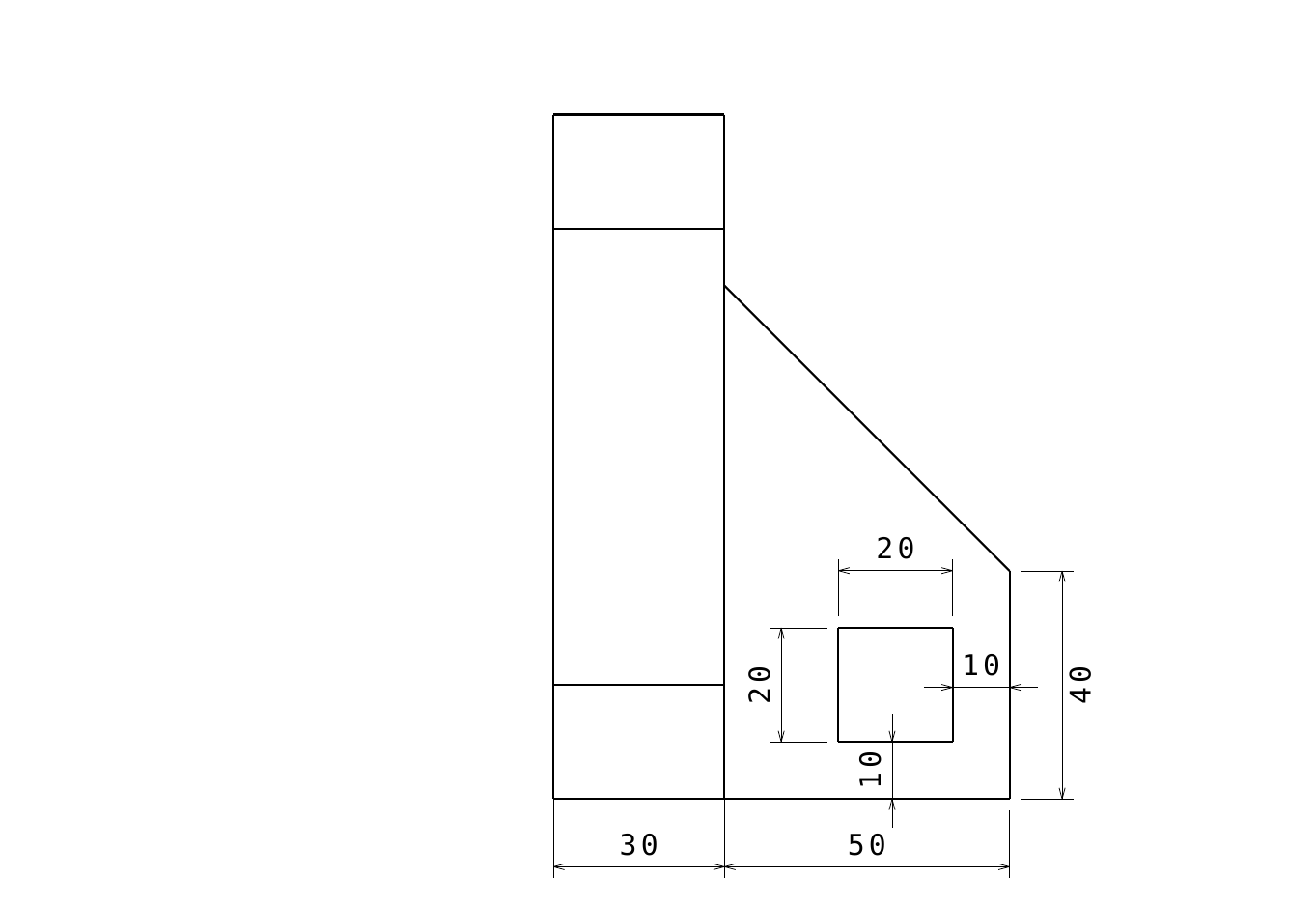}\label{fig:GE_bracket_tdb}}
	\caption{Engine bracket. Technical drawings.}
	\label{fig:GE_bracket_td}
\end{figure}

%
As a final example, we demonstrate the optimisation of a bracket lattice structure (loosely motivated by the GE jet engine bracket challenge) to demonstrate the application of the proposed approach to large problems.  We simplify the original engine bracket geometry while retaining its essential features, see Figures~\ref{fig:1_summary} and ~\ref{fig:bracketgroundside}. In particular, the four holes for fixing the horizontal base plate and the hole for applying the load to the vertical plate are omitted. The dimensions of the bracket are given in Figure~\ref{fig:GE_bracket_td}. The lattice consists of body-centred cubic (BCC) unit cells of size $5 \times 5$, leading to 82124 members and 12971 joints.

In the finite element model, all the joints on the reentrant surfaces of the base plate are fixed, and a unit load is applied to the 62 joints at the top left corner of the hole in the vertical plate.  As illustrated in Figure \ref{fig:bracketgroundside}, each unit load has a slope of $45\,\deg$ with respect to the horizontal base plate. The prescribed maximum volume is ${V_{\max}=15200}$, and the cross-sectional areas are constrained not to exceed ${A_{\max} = 1.0}$. The relative member densities $\tilde s_e (s_e)$ are penalised using the penalisation function introduced in Figure~\ref{fig:spline_plots}. Furthermore, the mean Young's modulus of the members is $\overline{E}=100$, and for the SPDE-defined covariance matrix the length-scale is $\ell = 10$, the smoothness parameter is $\nu=1.5$, and the standard deviation is $\sigma = 0.1\overline{E}$.

\begin{figure*}
	\centering
	\subfloat[$\alpha = 1$]{\includegraphics[width=0.33\linewidth]{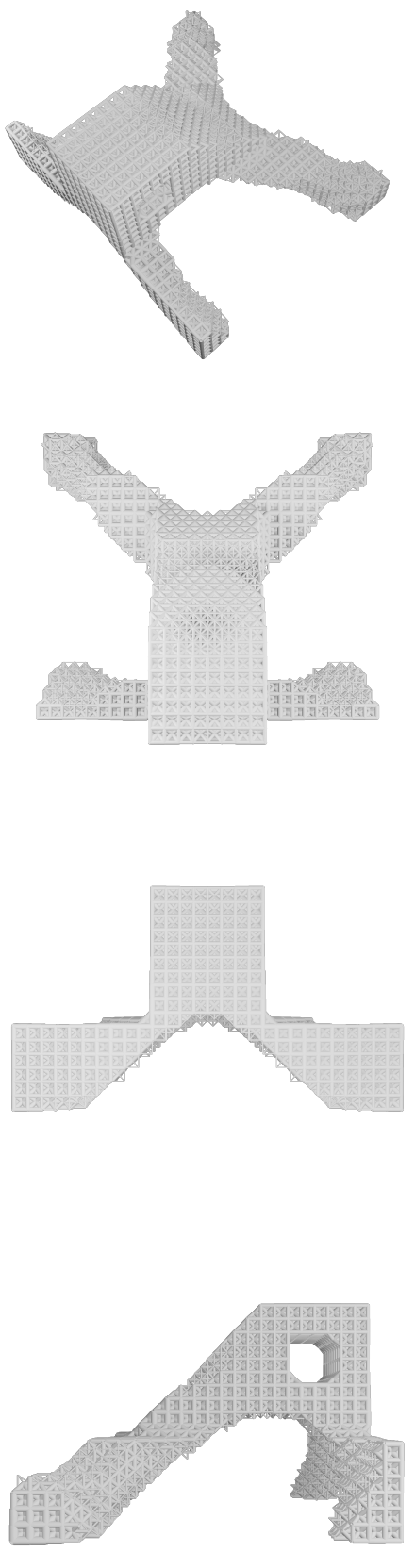}\label{fig:bracket_dso_all}}%
	\subfloat[$\alpha = 0.2$]{\includegraphics[width=0.33\linewidth]{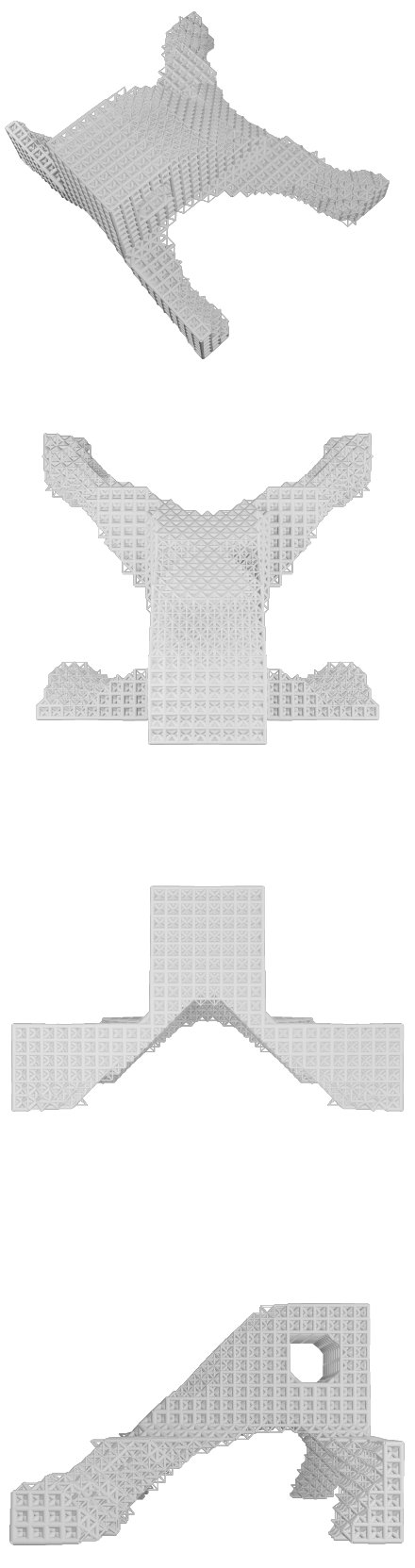}\label{fig:bracket_rso02_all}}%
	\subfloat[$\alpha = 0$]{\includegraphics[width=0.33\linewidth]{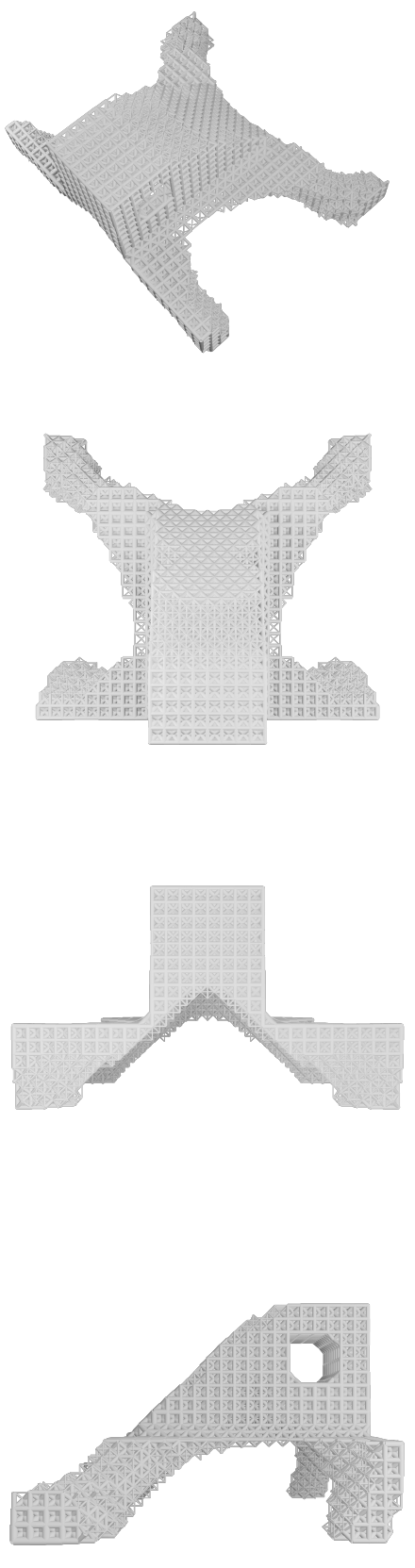}\label{fig:bracket_rso_all}}
	\caption{Engine bracket. Comparison of the optimised engine bracket for three different cost function weight parameters. For $\alpha = 1$ (left column) the cost function consists only of the mean and for $\alpha=0$ (right column) only of the standard deviation of the compliance. }
	\label{fig:bracket_compare}
\end{figure*}

Figure~\ref{fig:bracket_compare} shows the optimised lattice structures for three different cost function parameter values \mbox{$\alpha \in \{ 1, \, 0.2, \, 0 \}$}. As discussed, a decrease in the value of~$\alpha$  leads to increased consideration of the standard deviation of the compliance in the cost function. For the present structure, a decrease in~$\alpha$ results in a concentration of material around the top ends of the four legs of the bracket. For instance, reducing from~\mbox{$\alpha =1 $} to \mbox{$\alpha=0.2$} yields a reduction of the standard deviation  from 6.19 to 3.67 and an increase of the mean by 5\%. 

At last, in terms of computational efficiency, in each iteration step, it is sufficient to factorise the stiffness matrix of the truss structure only once. Subsequently,  the derivative of the displacement vector with respect to the member Young's moduli is computed using back substitution, cf.~\eqref{eq:displacement_derivative}. For each of the 82124  members, it is necessary to apply one back substitution which may represent a performance bottleneck in the case of large problems. However, each member's derivative of the displacement vector can be independently computed so that this step is easily parallelisable. \ADDED{In robust optimisation with $\alpha = 0.2$, the CPU run-time for one iteration is 18106.19\,s (1783.77\,s wall time), taking in total 127 iterations. For comparison, the run-time for a coarse lattice with 3226 members is 16.56\,s (3.44\,s wall time), and for a lattice with 10676 members it is 200.72\,s (38.53\,s wall time). The optimisation results for the two coarse meshes are not included in the paper. All experiments were run on an Intel\textregistered Xeon\textregistered Silver 4116 CPU @ 2.10\,GHz (12 cores, 24 logical processors).}

%
\section{Conclusions}
%
We introduced an efficient approach for topology optimisation of large lattice structures with spatially correlated random material properties. Although we considered only members with random Young's moduli, it is straightforward to adapt the proposed approach for other random fields, including the member cross-sectional areas, load distribution and joint positions, by following~\cite{guest2008structural}. We obtain the spatially correlated random fields by generalising the SPDE formulation of Mat\'ern fields to lattice structures. In the SPDE formulation, the random field is defined as the solution of a possibly fractional PDE with a second-order elliptic operator and a Gaussian white noise as the forcing. The properties of the resulting random field are governed by the parameters of the SPDE, i.e. the length-scale, the exponent and the standard deviation of the white noise. We generalise the SPDE formulation to lattices by replacing the stiffness and mass matrices of the finite element discretised SPDE with the respective matrices of the lattice structure. The resulting system matrix is the precision matrix of the Gaussian random field, which is usually sparse or can be approximated as sparse in the case of fractional SPDEs. We hypothesise that the SPDE models the unaccounted physical diffusion processes during the manufacturing of the lattice structure and, hence, is an effective approach to constructing physics-informed random fields. This interpretation of the SPDE formulation of random fields helps generalise them to non-stationary and anisotropic random fields.

For robust topology optimisation by optimising the cross-sectional areas of the members, we use a standard gradient-based approach and compute the required mean and standard deviation of the compliance and constraint functions using a first-order Taylor series expansion. The efficiency of the series expansion combined with the sparsity of the precision matrix of the random field allows us to optimise large lattice structures with up to eighty thousand members with random Young's moduli. These computations are only feasible because the expectation and standard deviation of the cost and the constraint functions and their gradients are computed using sparse matrix operations. 

We note that it is straightforward to apply the proposed robust optimisation approach to other types of lattice structures, including lattice-skin systems composed of a lattice core and a thin-shell skin~\citep{xiao2019interrogation, xiao2022infill} or structures consisting of beam members connected by joints that can transfer both forces and moments~\citep{yin2020CADcompatible}. Furthermore, we assumed in the present paper that the parameters of the SPDE are given. They will depend on the particulars of the manufacturing process, like the deposition rate and laser power in additive manufacturing by selective laser sintering~\citep{gibson2021additive}.  Relating the SPDE parameters to such directly observable manufacturing process parameters is of great practical relevance for designing more robust structures.

\section*{Acknowledgements}
I.B-Y. was supported by the European Union's Horizon 2020 research and innovation programme under the Marie Sklodowska-Curie Grant Agreement No. 101007815. A.O.Y. was supported by a fellowship of the Ministry of National Education of the Republic of Turkey. F. C. was supported by Wave 1 of The UKRI Strategic Priorities Fund under the EPSRC Grant EP/T001569/1 and The Alan Turing Institute.

\section*{Conflict of interest}
On behalf of all authors, the corresponding author states that there is no conflict of interest.

\section*{Replication of results}
The data that support the findings of this study are available from the corresponding author upon request.

\appendix
\ADDED{
%
\section{Finite element discretisation of the SPDE \label{appen:SPDEdiscret}}
%
In this appendix we briefly summarise the finite element discretisation of the recursive system of SPDEs~\eqref{eq:spdeRecurs} and the derivation of the probability density of the corresponding random field. 
See~\cite{lindgren2011explicit,koh2023stochastic} for a  more detailed discussion. 
}

\ADDED{
It is sufficient to detail the discretisation of the first SPDE~\eqref{eq:spdeRecurs1} in the recursion to discretise~\eqref{eq:spdeRecurs}. Its weak form reads: find~$r^{(1)} (\vec x) \in \set H^1 ( \Omega)$ such that 
\begin{equation}
	\begin{aligned}
	&\int_{\Omega} \left ( \kappa^2  r^{(1)} (\vec x ) v (\vec x)  +   \nabla r^{(1)} (\vec x) \cdot  \nabla v  (\vec x) \right)  \D \vec x \\
	&=  \frac{1}{\tau}  \int_{\Omega}  g(\vec x)v  (\vec x) \D \vec x \, \quad \forall \ v(\vec x)  \in \set H^1 (\Omega) \, ,
	\label{eq:spdeWeakForm}
	\end{aligned}
\end{equation}
where~$ \nabla$ is the gradient operator, $v(\vec x)$ is a test function and~$\set H^1 (\Omega)$ is a standard Sobolev space. The domain of~\eqref{eq:spdeRecurs} is all of~$\mathbb R^d$ which we approximate with a sufficiently large~$\Omega \subset \mathbb R^d$.
}

\ADDED{
To obtain the finite element approximation,~$r^{(1)} (\vec x)$ and~$v(\vec x)$ are discretised as
\begin{equation} \label{eq:discIntExp}
\begin{aligned}
r^{(1)} (\vec x) &\approx  \sum_{i} \phi_i (\vec x) r_{i}^{(1)} \, ,  \\
v  (\vec x) & \approx  \sum_{i} \phi_i (\vec x) v_i   \, , 
\end{aligned}
\end{equation}
where~$\phi_i (\vec x)$ are the (global) finite element basis functions, and~$r_{i}^{(1)}$ and $v_i$ are the respective nodal coefficients.  The so discretised weak form yields after quadrature 
\begin{equation} \label{eq:appPDErecursive1} 
\left ( \kappa^2 \vec M  + \vec A   \right )  {\vec r}^{(1)} = \frac{1}{\tau} {\vec g} \, , 
\end{equation}
where~$\vec M$ is the mass matrix,~$\vec A$ is the discretised Laplacian matrix  and~$\vec g$ is a random vector with the components
\begin{equation}
g_i = \int_{\Omega} g(\vec x) \phi_i(\vec x)  \D \vec x  \, .
\end{equation}
Here,~$g(\vec x)$ is a Gaussian white noise process with unit variance as given in~\eqref{eq:GPwhiteNoise}. After defining the matrix~$\vec L =  \kappa^2 \vec M  + \vec A $, we can write~\eqref{eq:appPDErecursive1} as
\begin{equation} \label{eq:appPDErecursive1compact} 
	\vec L  {\vec r}^{(1)} = \frac{1}{\tau} {\vec g} \, . 
\end{equation}
The discretisation of the subsequent equation in the recursion~\eqref{eq:spdeRecurs2} proceeds similarly, yielding
\begin{equation} \label{eq:spdeFErecurs2b}
   \vec L \vec r^{(k)}  = \M \vec r^{(k-1)} \, .
\end{equation}
}

\ADDED{
Next, we consider the probability density of the random solution~$\vec r$. The random source vector~$\vec g$ in~\eqref{eq:appPDErecursive1} is a zero-mean Gaussian random vector given that~$g (\vec x)$ is a zero-mean Gaussian white noise process with unit variance. Hence, the components of covariance matrix of~$\vec g$ are given by
\begin{align}
\begin{split}
&\cov \left (  g_i , g_j \right )  = \\
&= \expect \left [  \int_{ \Omega}  \int_{ \Omega}   g(\vec x) \phi_i(\vec x)     g(\vec x') \phi_j(\vec x')  \D \vec x  \D \vec x'   \right ] \\
&=  \int_{\Omega} \int_{ \Omega}  \phi_i(\vec x)  \expect \left [   g(\vec x)   g(\vec x') \right ] \phi_j(\vec x') \D \vec x \D \vec x'  \\
&=  \int_{ \Omega} \int_{\Omega}  \phi_i(\vec x)  \delta (\vec x - \vec x') \phi_j(\vec x') \D \vec x \D \vec x'   \\
&=  \int_{ \Omega} \phi_i(\vec x) \phi_j(\vec x) \D \vec x \, . 
\end{split}	
\end{align}
Thus,~$\vec g$ has the multivariate Gaussian density 
\begin{equation}\label{eq:MassDens}
p( \vec g ) = \set N (\vec 0, \vec M) \, .
\end{equation}
According to the linear transformation property of Gaussian densities, see e.g. \cite{williams2006gaussian}, and~\eqref{eq:appPDErecursive1compact} the random vector~$\vec r^{(1)}$ has the density
\begin{equation}
p (  \vec r^{(1)} ) =  \set N \left (\vec 0, \, \vec C_{r}^{(1)}  \right ) =  \set N \left (\vec 0, \,  \frac{1}{\tau^2}\vec L^{-1} \vec M \vec L^{-\trans} \right )   \, .
\end{equation}
The respective precision matrix is
\begin{equation}
	\vec Q_{r}^{(1)} = \left ( \vec C_{r}^{(1)} \right )^{-1}  = \tau^2 \vec L^\trans \vec M^{-1} \vec L \, .
\end{equation}
In the same way, we obtain using the recursion~\eqref{eq:spdeFErecurs2b}  and the linear transformation property of Gaussian densities the recursion equation for the precision matrix 
\begin{equation}
    \vec Q_{r}^{(k)} = \vec L^{\trans} \vec M^{-\trans}\vec Q_{r}^{(k-1)}\vec M^{-1}\vec L\, . 
\end{equation}
The two preceding equations can be combined to obtain~$\vec Q_{r} \equiv  \vec Q_{r}^{(\beta)}$ given in~\eqref{eq:precision}.  In our computations the lumped mass matrix is sparse so that the precision matrix~$\vec Q_{r}$ is sparse.
}

\bibliography{robust-opt}

\end{document}